\documentclass[a4paper,10pt,fleqn]{article}%

\usepackage{blindtext}

\usepackage[numbers,square]{natbib}%
\usepackage[bottom=3cm,top=2cm,left=2cm,right=2cm]{geometry}%
\usepackage{framed,url,setspace,graphicx,bm,xcolor,siunitx,nicefrac}
\usepackage{amsmath,amssymb,amsthm}%
\usepackage{stmaryrd,textcomp,colortbl,ifthen}%

\usepackage{placeins}

\usepackage[utf8]{inputenc}

\newcommand{\thisisrevision}[1]{#1}
\newcommand{\thisfigurepath}{./figures/}

\usepackage{notation}
\usepackage{mathoperators}

\usepackage{pgfplots}%
\pgfplotsset{width=7cm, compat=1.6, x tick label style={/pgf/number format/.cd,fixed,set thousands separator={}},%
y tick label style={/pgf/number format/.cd,fixed,set thousands separator={}}}%

\newcommand{\reffig}[1]{Figure~\ref{fig:#1}}%
\newcommand{\refFig}[1]{Figure~\ref{fig:#1}}%
\newcommand{\reftab}[1]{Table~\ref{tab:#1}}%
\newcommand{\refEqn}[1]{Equation~\ref{eqn:#1}}%
\newcommand{\refeqs}[1]{Eqs.~\ref{eqn:#1}}%
\newcommand{\refeqno}[1]{\ref{eqn:#1}}%
\makeatletter
\newcommand{\refeqn}{\@ifstar{\@ifstar\@@refeqn\@refeqn}{\@@@refeqn}}%
\def\@refeqn#1{(\ref{eqn:#1})}%
\def\@@refeqn#1{eqs.~(\ref{eqn:#1})}%
\def\@@@refeqn#1{eq.~(\ref{eqn:#1})}%
\makeatother

\newcommand{\refeqnsafe}[1]{eq.~(\ref{eqn:#1})}%

\colorlet{notecolor}{black!12}%
\colorlet{remarkcolor}{black!12}%
\newtheorem{prototheorem}{Note}[section]
\newenvironment{note}%
{\colorlet{shadecolor}{notecolor}\begin{shaded}\begin{prototheorem}}%
{\end{prototheorem}\end{shaded}}%
\newcommand{\refnote}[1]{note~\ref{note:#1}}%

\newtheorem{remarktheorem}{Remark}[section]%
\newenvironment{remark}%
{\colorlet{shadecolor}{remarkcolor}\begin{shaded}\begin{remarktheorem}}%
{\end{remarktheorem}\end{shaded}}%
\usepackage{reconstructnotation}%

\usepackage{hyperref}

\def\updatefigure#1{}

\updatefigure{380}%
\updatefigure{381}%
\definecolor{newcolor}{rgb}{.8,.349,.1}

\def\worktitle{Efficient sequential PLIC interface positioning for enhanced performance of the three-phase VoF Method}%

\providecommand{\keywords}[1]{\textbf{\textit{Keywords---}} #1}

\allowdisplaybreaks

\makeatletter
\def\underbracex#1#2{\mathop{\vtop{\m@th\ialign{##\crcr
   $\hfil\displaystyle{#2}\hfil$\crcr
   \noalign{\kern3\p@\nointerlineskip}%
   #1\crcr\noalign{\kern3\p@}}}}\limits}

\def\underbracea{\underbracex\upbracefilla}

\def\upbracefilla{$\m@th \setbox\z@\hbox{$\braceld$}%
  \bracelu\leaders\vrule \@height\ht\z@ \@depth\z@\hfill 
\kern\p@\vrule \@width\p@\kern\p@\vrule \@width\p@\kern\p@\vrule \@width\p@
$}

\def\upbracefillb{$\m@th \setbox\z@\hbox{$\braceld$}%
\vrule \@width\p@\kern\p@\vrule \@width\p@\kern\p@\vrule \@width\p@\kern\p@
 \leaders\vrule \@height\ht\z@ \@depth\z@\hfill\bracerd
  \braceld\leaders\vrule \@height\ht\z@ \@depth\z@\hfill
\kern\p@\vrule \@width\p@\kern\p@\vrule \@width\p@\kern\p@\vrule \@width\p@
$}

\def\upbracefilld{$\m@th \setbox\z@\hbox{$\braceld$}%
\vrule \@width\p@\kern\p@\vrule \@width\p@\kern\p@\vrule \@width\p@\kern\p@
 \leaders\vrule \@height\ht\z@ \@depth\z@\hfill\braceru$}

\def\underbracebd{\underbracex\upbracefillbd}
\def\upbracefillbd{$\m@th \setbox\z@\hbox{$\braceld$}%
\vrule \@width\p@\kern\p@\vrule \@width\p@\kern\p@\vrule \@width\p@\kern\p@
\bracerd\braceld
 \leaders\vrule \@height\ht\z@ \@depth\z@\hfill\braceru$}
 
\makeatother

\begin{document}
\newcommand{\refapp}[1]{appendix~\ref{app:#1}}%
\title{\worktitle}%
\author{%
Johannes Kromer\textsuperscript{1}\href{https://orcid.org/0000-0002-6147-0159}{\includegraphics[height=10pt]{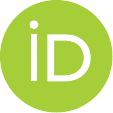}}%
, %
Johanna Potyka\textsuperscript{2,$\dagger$}\href{https://orcid.org/0000-0003-1310-4434}{\includegraphics[height=10pt]{\thisfigurepath orcid_logo}}%
, %
Kathrin Schulte\textsuperscript{2}\href{https://orcid.org/0000-0001-8650-5840}{\includegraphics[height=10pt]{\thisfigurepath orcid_logo}}%
\ and %
Dieter Bothe\textsuperscript{1}\href{https://orcid.org/0000-0003-1691-8257}{\includegraphics[height=10pt]{\thisfigurepath orcid_logo}}%
}%
\date{}%
\maketitle
\begin{center}
\small\em%
\textsuperscript{1}%
Mathematical Modeling and Analysis, Technische Universit\"at Darmstadt\\ Alarich-Weiss-Strasse 10, 64287 Darmstadt, Germany\\%
\textsuperscript{2}%
Institut f\"ur Thermodynamik der Luft- und Raumfahrt, Universit\"at Stuttgart\\ Pfaffenwaldring 31, 70569 Stuttgart, Germany\\%
\textsuperscript{$\dagger$}Email for correspondence: %
\href{mailto:johanna.potyka@itlr.uni-stuttgart.de?subject=Efficient\%20three-phase\%20PLIC}{johanna.potyka@itlr.uni-stuttgart.de}%
\end{center}

\let\paragraphold\paragraph%
\renewcommand{\paragraph}[1]{\paragraphold{#1.}}%

\begin{abstract}
%
%
%
This paper presents an efficient algorithm for the sequential positioning, also called nested dissection, of two planes in an arbitrary polyhedron. Two planar interfaces are positioned such that the first plane truncates a given volume from this arbitrary polyhedron and the next plane truncates a second given volume from the residual polyhedron. %
This is a relevant task in the numerical simulation of three-phase flows when resorting to the geometric Volume-of-Fluid (VoF) method~\cite{HIRT1981} with a Piecewise Linear Interface Calculation (PLIC)~\cite{Youngs1982}. %
An efficient algorithm for this task significantly speeds up the three-phase PLIC algorithm. %
The present study describes a method based on a recursive application of the \textsc{Gaussian} divergence theorem, where the fact that the truncated polyhedron shares multiple faces with the original polyhedron can be exploited to reduce the computational effort. %
A careful choice of the coordinate system origin for the volume computation allows for successive positioning of two planes without reestablishing polyhedron connectivity. %
Combined with a highly efficient root finding, this results in a significant performance gain in the reconstruction of the three-phase interface configurations. %
\\
The performance of the new method is assessed in a series of carefully designed numerical experiments. %
Compared to a conventional decomposition-based approach, the number of iterations and, thus, of the required truncations was reduced by up to an order of magnitude. %
The PLIC positioning run-time was reduced by about $90\%$ in our reference implementation. %
Integrated into the multi-phase flow solver Free Surface 3D (FS3D), an overall performance gain of about $20\%$ was achieved. %
The reference implementation of the efficient sequential PLIC positioning is available as a Fortran module at \href{https://doi.org/10.18419/darus-2488}{https://doi.org/10.18419/darus-2488}. %
Allowing for simple integration into existing numerical schemes, the proposed algorithm is self-contained, requiring no external decomposition libraries. %
%

\end{abstract}

%
%
%
\section{Introduction}\label{sec:introduction}%
Multiphase flows with three phases play an important role in nature, science, and technology, ranging from processes in clouds where droplets interact with already frozen particles to groundwater predictions where water and air interact with the solid phase in porous media. %
As the examples indicate, not only two but often three phases interact. %
Predicting technical processes for optimization and forecasting natural phenomena is required in multiple engineering tasks, and numerical simulations are one of the main tools to achieve this. %
High efficiency of the involved simulation tools is not only advantageous in increasing the speed of the simulation, but also required to reach an efficient use of limited resources. %
A small selection of three-phase flow simulations shows the wide range of applicability and relevance: \citet{PATEL2017} simulate an idealized pore structure formed by spheres which is filled with oil and flooded with water. %
\citet{BAGGIO2019} on the other hand studied the impact of a droplet on a structured surface, while \citet{REITZLE2017} investigated the solidification process of water, where three-phase contact lines are present shortly before the droplet is fully frozen. %
The encapsulation process after the head-on collision of two droplets of immiscible liquids was simulated by \citet{Li2015}, where the three-phase contact line is formed by two liquids instead of solid and liquid in an ambient phase. %
Due to its high practical relevance, we thus specifically enhance the \textbf{three phase} case. %

The Volume of Fluid method (VoF) by \citet{HIRT1981} is widely employed in numerical simulations investigating multi-phase flow. %
The advection of the phase indicator field, discretized as the volume fractions, maintains the identification of the phases throughout the simulations. %
An accurate advection of the volume fractions requires additional methods to maintain a sharp interface: %
In geometric VoF methods, cf.\ \citet{MARIC2020} for a recent review, the \textbf{Piecewise Linear Interface Calculation (PLIC)} often performs this task, by locally, i.e., in every cell containing a piece of interface, placing a plane inside the mesh cell such that the local volume fractions are matched, while the orientation of the plane is chosen from a suitable normal calculation. %
In three-phase flows, each individual phase interface can be approximated by positioning of a local plane, hence up to three PLIC-planes can be present inside a single mesh cell. %
Such a \textbf{three-phase PLIC reconstruction} avoids numerical smearing as the subsequent advection can be performed via geometrical flux computations. %

\begin{figure}[!tb]%
\centering
\includegraphics[width=0.5\textwidth]{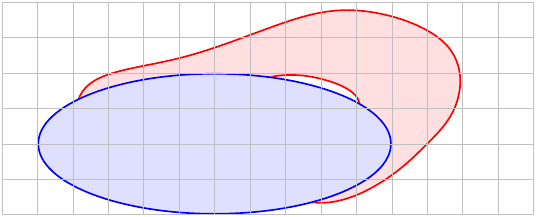}%
\caption{Interfaces intersecting at triple lines (here: triple points): Resulting multi-material phase distributions with wetted and non-wetted configurations may display thin gas or liquid films.}%
\label{fig:topology_illustration}%
\end{figure}

\begin{figure}
\null\hfill%
\includegraphics[page=4,scale=0.8]{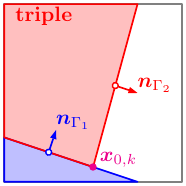}%
\hfill
\includegraphics[page=5,scale=0.8]{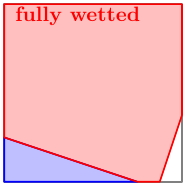}%
\hfill
\includegraphics[page=3,scale=0.8]{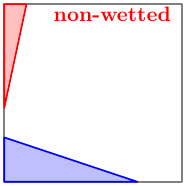}%
\hfill\null%
\caption{Different configurations in three-phase cells: Configurations with triple lines, fully wetted and non-wetted are possible. (For simplicity, the two dimensional analoga are displayed.) \vspace{\baselineskip}}%
\label{fig:threephasecells}%
\end{figure}

As illustrated in \reffig{topology_illustration}, a variety of different phase topologies occur in the discretized problem of three-phase interaction: Contact lines and thin films of liquid or gas are possible within a three-phase cell. %
Thus, the PLIC method has to be able to reconstruct the different three-phase cell configurations illustrated in \reffig{threephasecells}. %
At the intersection of three interfaces, the so called triple line, the involved interfaces either enclose a prescribed, e.g.,~\citet{PATEL2017}, or an iteratively determined, e.g.,~\citet{JCP_2016_atdv}, angle, usually referred to as the contact angle in solid-liquid-gas interaction. %
Besides a triple configuration, a three-phase cell may also admit a non-wetted or fully-wetted configuration, which can also contain a film of liquid or gas thinner than the cell height. %

\begin{figure}[!tb]%
\null\hfill%
\includegraphics[scale=0.96,page=3]{./figures/literature_onionskin}%
\hfill%
\includegraphics[scale=0.96,page=1]{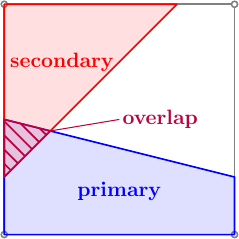}%
\hfill\null
\\%
\null\hfill
\includegraphics[scale=0.96,page=2]{./figures/literature_onionskin}%
\hfill%
\includegraphics[scale=0.96,page=5]{./figures/literature_onionskin}%
\hfill\null%
\caption{Different algorithms for multi-material reconstruction: In the top row, configurations are shown, which do not consider the presence of a triple line: The onion-skin method with a layer structure (top left) and unphysical domain overlap which can be present both in the onion-skin method or a completely independent reconstruction of the phases (top right). In the second row methods are depicted which consider the presence of a triple line: nested dissection (bottom left) and power-diagram with generators $\vg_i$ and weights $\omega_i$ (bottom right). Note that, with corresponding normals and volume fractions, the nested dissection redistributes the domain overlap produced by the onion-layer scheme (purple, hatched).}%
\label{fig:algorithm_concepts}%
\end{figure}

The three-phase PLIC interface reconstruction consists of three tasks: \textbf{material ordering}, \textbf{interface normals orientation}, and \textbf{interfaces' positioning}. %
The interfaces are placed such that the prescribed volumes (in 3D) or areas (in 2D) given by the previous VoF advection step are enclosed. %
\refFig{algorithm_concepts} provides an overview of the conceptual differences for three-phase PLIC which are discussed in some more detail in the literature overview in the following. %
The methods are ordered from left to right with respect to increasing accuracy of the method in representing the true physical situation, but also increasing computational complexity and effort: 
The \textbf{onion-skin} (e.g.\ \cite{Youngs1982, Benson2002}) approach consists of a sequential positioning of all planes in the polyhedral cell: %
With $N$ pairs of volume fractions and normals $\set{(\refvof,\plicnormal)}_{n=1}^{N}$, the $n$-th plane is positioned such that it encloses the sum of all volumes associated with the previous $(n-1)$ planes and the additional $n$-th volume fraction. %
This method ignores the location of the previously positioned planes. %
The simplicity in concept and implementation, however, also bears a major disadvantage: The method can represent a layered configuration correctly, but some combinations of input data may lead to an unphysical overlap of the reconstructed phase domains, cf.\ the second panel in \refFig{algorithm_concepts}. %
The \textbf{sequential or nested dissection} (e.g.\ \cite{IJNMF_2007_apla,MCM_2008_anmf,JCP_2016_atdv}) overcomes this unphysical artifact by positioning the $n$-th plane in a polyhedron truncated by its $n-1$ predecessors; see the illustration in the third panel of \reffig{algorithm_concepts}. %
In the onion skin and in the sequential PLIC method the number of planes representing the two (no triple line) or three (with a triple line) interfaces in three-phase cells is reduced to two planes: the solid's or other primary phase's interfaces towards both other phases are assumed to have no kinks which greatly reduces the algorithm's complexity. Being a reasonably compromise of accuracy and efficiency, the sequential PLIC method is chosen for the present work. %
The \textbf{power-diagram} (cf.\ \cite{JCP_2009_asoa}) method introduces an additional degree of freedom to the orientation and positioning task by allowing for kinks in the respective interfaces; see the rightmost panel in \reffig{algorithm_concepts}. It is the most accurate, but in three dimensions an infeasibly time-consuming method.%

Because of the inherent difference in mathematical and implementation complexity, the following literature overview of three-phase PLIC is divided based on the spatial dimensions. %
%
%
\subsubsection*{Two spatial dimensions}%
%
%
%
\citet{Youngs1982} described the ``onion-skin" method for the treatment of multiple phases in a Cartesian grid. %
The two-phase PLIC positioning like described in detail by \citet{Benson2002} for cuboid cells is sufficient to fulfil the positioning tasks. %
The major advantage is the simplicity, as the same positioning as for two phases is used. %
The major drawback of this method is that overlaps are possible, thus during the subsequent advection either mass conservation is not ensured, if the flux is limited to the volume of the advected box, or overfull cells may be produced, i.e. cells which contain more volume of the advected phases than the cell's volume.
 %

\citet{IJNMF_2007_apla}~introduce ``an approach to reconstructing a three-material cell that accurately estimates interface normals, determines the position of a triple point if it exists, and exactly conserves mass'' for \textsc{Cartesian} meshes. Employing a nested-dissection type approach, they require a pre-determined sequence of reconstruction, where the first interface normal is computed by the two-phase PLIC method of~\citet{Youngs1982}. %
The second normal is reconstructed such that the volume fractions predicted by the corresponding plane best approximate (in a least-squares sense) those in the neighboring mesh cells. %

%
\citet{MCM_2008_anmf} employ a minimization procedure for the reconstruction of a triple point based on the neighboring cell normals. %
In both of the aforementioned studies the employed positioning is not described further beyond choosing the position such that it encloses the given area (analogous to the volume in 3D). %

%
\citet{JCP_2009_asoa}~propose a second-order accurate (with respect to the $L_2$-norm) material-order-independent interface reconstruction technique, based on the work of~\citet{IJNMF_2008_moii}, where no material ordering has to be pre-determined. %
They employ a weighted \textsc{Voronoi} diagram (called \textit{power diagram}) to simultaneously approximate the phase-centroids of all materials, which are then used to determine the respective normals on equidistant \textsc{Cartesian} meshes; cf.\ the rightmost panel of~\reffig{algorithm_concepts}. %
By iteratively adjusting the weights of the so-called point generators associated to the materials in a cell, the interfaces are positioned to match the respective volume fractions. %
First, those point generators $\vg_i$ are computed by solving a least-squares problem in a rectangular domain. %
With those generators $\vg_i$ fixed in space, the weights $\omega_i$ are adjusted such that the induced volume fractions match the prescribed ones. %
However, to the best of our knowledge, all published methods of this class were presented for two space dimensions, and the complexity can be expected to increase strongly in three spatial dimensions \cite{JCP_2009_asoa}. %
\citet{JCP_2009_asoa} also found that ``for cells with high aspect ratios, this can [...] give poor results for multi-material cells''. %
\citet{IJNMF_2012_eoeo} conduct a series of numerical tests for the aforementioned algorithm. %

%
\citet{CGF_2008_dmmi}~reconstruct the interfaces in \textsc{Cartesian} cells by a combination of subdivision (between 125 and 1000 per cell) and probabilistic reordering by simulated annealing, based on the volume fraction data of the neighboring cells. %
While, for the cases they investigate, the reconstructed interfaces are qualitatively superior to those obtained by PLIC, the huge computational effort induced by the subdivision renders their algorithm impractical, especially for an extension to three dimensions. %
%
%
\subsubsection*{Three spatial dimensions}%
%
%
The extension of the "onion-skin´´ method to multiple phases in three dimensions is analogous to the two dimensional case. %
\citet{Benson2002} described the two-phase positioning for two-phases not only in two, but also in three dimensions, thus also the positioning for a three-dimensional ``onion-skin" interface reconstruction. %
As already mentioned above, the ``onion-skin" approach leads to overlaps in cases with a triple line, leading to poor mass conservation during a subsequent advection. %
%
%
Various authors, e.g.,~\citet{PATEL2017, Washino2010}, employ a sequential PLIC approach, cf.\ the third panel in \reffig{algorithm_concepts}, where the primary solid phase interface is reconstructed independently of the other two phases by Youngs~\cite{Youngs1982} two-phase PLIC algorithm. %
In three-phase cells, the second interface dividing the second and third phase is reconstructed subsequently by imposing a static- or dynamic contact angle from a dynamic contact angle correlation in all three phase cells.
This is achieved by imposing a user defined or pre-determined contact angle between the primary and secondary PLIC planes. %
In a first step, the normal is obtained from Youngs's method and, in a second step, the orientation towards the primary plane is adjusted to match the contact angle. %
Note that the overall accuracy of this methods also strongly depends on the quality of the contact angle correlation for the simulated problem. %
On the one hand, this method is not able to reconstruct thin films, as the contact angle is enforced in all three-phase cells. On the other hand, the normals of the primary and secondary planes are computed explicitly and thus this algorithm is relatively fast compared to iterative methods. %
\citet{PATEL2017}~employ ``transformations of the coordinate system [...] to reduce the number of possible configurations of a PLIC surface from 64 to 5 generic configurations'' and refer to the method by~\citet{VANSINT2005} for the positioning which is restricted to \textsc{Cartesian} applications. This method is only described for two phases, thus it remains unclear, if an overlap of the PLIC interfaces is allowed in the positioning by~\citet{PATEL2017} or if this is somehow incorporated into the mentioned transformation of the cuboid cell.

%
\citet{JCP_2016_atdv}~introduce a comprehensive algorithm for the sequential three-phase PLIC interface reconstruction, covering both normal reconstruction and positioning of the PLIC interfaces from given volume fractions. %
The latter part of their algorithm resorts to a combination of polyhedron decomposition and usage of \textsc{Brent}'s root-finding method to match the target volumes. %
\citet{JCP_2016_atdv}~report better than first-order convergence for the computed normals (with respect to the $L_2$-norm). %
Information on the convergence of the positioning algorithm based on the clipping and capping algorithm by \citet{Stephenson1975} is not reported. %
Also the original description of the clipping and capping does not provide information on the number of iterations or other performance data. %

\begin{figure*}[!tb]%
\null\hfill%
\includegraphics[page=1]{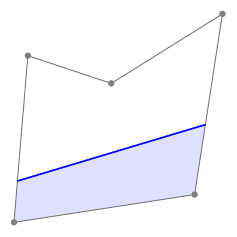}%
\hfill%
\includegraphics[page=2]{\thisfigurepath literature_shashkov}%
\hfill%
\includegraphics[page=3]{\thisfigurepath literature_shashkov}%
\hfill%
\includegraphics[page=4]{\thisfigurepath literature_shashkov}%
\hfill\null%
\caption{The nested (sequential) dissection algorithm of \protect\citet{JCP_2007_mmir} employs tetrahedron decomposition in three dimensions, the triangular decomposition in two dimensions is visualized.}%
\label{fig:illustration_shaskov}%
\end{figure*}

%
\citet{JCP_2007_mmir} introduce a recursive algorithm for sequential reconstruction based on tetrahedron decomposition, also called nested dissection, which is applicable to convex and non-convex polyhedra. %
However, for the latter to be accessible for their tetrahedron decomposition scheme, the weighted sum of the vertices needs to be a star-point%
\footnote{A star-point $\vx$ of a polyhedron is any point $\vx\in\polyhedron*$ such that $\vx+\lambda\brackets{\vy-\vx}\in\polyhedron*\,\forall(\lambda,\vy)\in[0,1]\times\polyhedron*$, i.e., the line connecting $\vx$ with any other point $\vy\in\polyhedron*$ must be entirely contained in the polyhedron $\polyhedron*$. %
The weighted sum of the vertices can be assumed to be a star-point for weakly non-convex polyhedra.}%
of the polyhedral cell. %
The computational effort of this strategy can grow very quickly; cf.\ \reffig{illustration_shaskov}. %

\vspace{\baselineskip}

To sum up, quite a number of different approaches for the normal computations in two- and three dimensions already exist. %
A few positioning algorithms for three-phase PLIC are also available, which all employ a tetrahedron decomposition for the volume computation and re-establishing connectivity after the first truncation. %
Avoiding both, the volumes' tetrahedron decomposition as well as the connectivity computations, can significantly speed-up implementations of the three-phase PLIC method. 

\subsubsection*{Objectives and Strategy}
An efficient, robust and accurate interface positioning method, applicable for three phases and -ideally- suitable for any three-dimensional polyhedron shape, seems to be lacking from literature. %
Therefore, the aim of this work is the development of an enhanced sequential PLIC positioning method, as the literature review identifies this concept as a compromise of accuracy and speed of the algorithm. %

We assume the material order as given: The primary interface is reconstructed first in each cell. %
We also assume the interface normals as given, as there is a variety of suitable algorithms available. %
In the prototypical application of the new positioning algorithm within the in-house multi-phase flow solver Free Surface 3D (FS3D) presented below, the normals are computed with an iterative algorithm. This allows the full spectrum of possible topological phase configurations shown in \reffig{threephasecells}. %
The algorithm chosen for the normals' computation in the reference application in the present work is similar to the iterative three-phase PLIC normal orientation method by \citet{JCP_2016_atdv}. %
Such an iterative algorithm for the secondary plane's orientation requires multiple computations of the interfaces' positions. %
Each orientation during the root-finding of the optimal orientation also requires the computation of a position matching the correct volume fraction. %
Thus, a more efficient algorithm for the three-phase PLIC positioning, which is also an iterative process, becomes even more advantageous for the speed-up of the simulation tool if employed multiple times within such an iterative algorithm. %
Consequently, this work focuses on the enhancement of the PLIC positioning in three-phase cells. %

The accurate and efficient positioning of the PLIC planes poses a non-trivial task: %
Firstly, the truncation of a polyhedron cell results in different sub-polyhedra; secondly, each volume truncated by the PLIC plane must match the prescribed target volume. %
When moving the plane through the cell along the given normal, different volumes are enclosed by the plane. %
We seek for the unique position of the plane, for which a given target volume, prescribed by the result of the previous VoF advection step, is matched. %
Efficiently reaching a match of the target volume up to a given volume tolerance requires a fast root-finding algorithm. 
The choice of the volume tolerance as the termination criterion for this root-finding is discussed further in Note~1.1. %
If more than one target volume has to be matched by more than one plane, the positioning task becomes considerably more demanding. %
The present work introduces an efficient iterative algorithm for the sequential positioning of the two planes representing up to three interfaces in three-phase cells. %

\begin{note}[Volume- vs.\ position tolerance]%
\label{note:volumepositiontol}
Volume conservation is a major advantage of finite volume schemes, thus choosing a tolerance for the enclosed volume instead of a position tolerance ensures a given precision of the volume conservation of the interfaces for the advection. %
The volume tolerance cannot be translated into a single position tolerance of the planes. %
Each deviation of a position $\delta t$ is coupled to a deviation of the volume fraction $\delta \alpha$ by the immersed area, which in turn depends on the orientation of an enclosed volume below the PLIC planes (see illustration). %
This does not allow a comparison of the exact position between different positioning approaches, since the proportionality factor, i.e., $|\polyhedron*\cap\Gamma|$ may change between cells. %

but the accuracy of the method is given by the achievable volume tolerance. %
In our evaluations of the positions' deviation from the previous method's results, both with the volume fraction tolerance $10^{-14}$, differences of the positions were below $10^{-12}$ for cells of size $1^3$. %
The overall PLIC accuracy depends primarily on the chosen normal computation. \\%
\null\hfill
\includegraphics{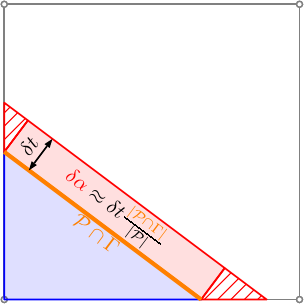}%
\hfill\null
\end{note}

The positioning problem at hand admits a strictly local character: Once the interface normals in a cell are approximated, the positioning problems in each cell containing two- or three phases are independent. %
However, a positioning scheme usually is embedded into a flow solver, thus the two- and three-phase positioning algorithm are typically used in parallel computations. An uneven distribution of single-, two- and three-phase cell types is usually present among the compute cores when a domain decomposition is employed, the parallelization strategy common in computational fluid dynamics. %
This implies that the slower execution time of a three-phase positioning algorithm is a bottleneck in the performance of the positioning in all cells and inducing load imbalances. %
A load rebalancing, i.e., a redistribution of the cells leading to a new domain decomposition, specific to suit the three-phase interface reconstruction would be costly in other parts of the flow solver. %
The multi-grid pressure solver employed to solve the pressure \textsc{Poisson} equation usually consumes the largest amount of the runtime in typical incompressible flow solvers. The multi-grid solver is most efficient with an evenly sized and cuboid-shaped domain per computer core. %
The most efficient approach is thus choosing the domain decomposition to suit the pressure solver and speed-up the three-phase positioning algorithm by a more efficient scheme like the one presented here. %
The integration into a flow solver avoiding the described load balance problem was a main goal of the present work besides the optimization of the stand-alone reconstruction method. %

In order to allow for a performance analysis of the new positioning algorithm not only as a stand-alone implementation, the latter is embedded into the in-house code FS3D \cite{EISENSCHMIDT2016}. %
FS3D is a massively parallelized multiphase flow solver employing \textsc{Cartesian} grids. %
While it is possible in a cuboid cell of such a \textsc{Cartesian} grid to compute the first positioning with two-phase PLIC, c.f. appendix \ref{app:cube_explicit_volume_inverse}, this leads to disadvantages arising for the subsequent second positioning: %
The second positioning has to be performed in the cut cell, the remaining polyhedron. %
In one form or another, the methods discussed in the literature review employ extraction of the topological connectivity of the cut cell and a decomposition of the volumes into tetrahedrons for the volume computation. %
A first idea to overcome at least the costly tetrahedron decomposition is performing the second positioning in the cut cell with the efficient positioning method for one plane in an arbitrary polyhedron by \citet{JCP_2021_fbip}. %
While this leads to a satisfying solution to the sequential positioning problem in terms of accuracy, and the performance is already increased compared to a decomposition based approach, the performance can be further enhanced: %
An improved sequential PLIC positioning presented in this study is possible, if some further considerations are taken into account. %

The improved algorithm does not perform two independent successive positioning steps, but makes use of the shared geometrical information from the first- and second cut. %
The fact that the truncated polyhedron shares multiple faces with the original polyhedron can be exploited to reduce the computational effort. %
Indeed, a careful choice of the coordinate system origin for the volume computation allows to compute the volume of the successive cuts without re-establishing polyhedron connectivity. %
With such a well-chosen origin, the interfaces formed by the intersections do not contribute to the volume computation and thus no connectivity has to be reestablished. %
This provides an additional performance gain compared to a successive application of two different efficient two-phase algorithms for cuboids and polyhedrons, which was, as elaborated on in the previous paragraph, the first idea to speed-up the three-phase PLIC algorithm in the reference application. %

Furthermore, the resulting algorithm is not limited to cuboid cells as the concept has no limitations concerning the shape of the polyhedron. %
The algorithm is applicable to any polyhedron cell geometry represented by a list of faces, including non-convex polyhedra. %
This widens the applicability of the efficient sequential three-phase PLIC positioning method to any geometric VoF framework with polyhedral cells simulating three-phase problems. This increases the relevance for applications beyond \textsc{Cartesian} solvers such as FS3D. %
The reference implementation of the efficient sequential three-phase PLIC positioning is made available as a Fortran module (\href{https://doi.org/10.18419/darus-2488}{https://doi.org/10.18419/darus-2488}), directly applicable to solvers with cuboid cells and easily adaptable to the needs of other codes. %
The volume computation with arbitrary polyhedral cells was tested extensively by Kromer and Bothe~\cite{JCP_2021_fbip}, which enforces this claim. %
Our implementation utilizes a list of faces for the reconstruction in three-phase cells, thus the cuboid can be replaced by any kind of polyhedron (including non-convex polyhedra). %

Already for the uniform Cartesian mesh used in FS3D, a significant performance gain of about 90\% for the three-phase reconstruction was achieved resulting in about 20\%  overall performance gain for the whole application. %
Solvers with more complex meshes composed of arbitrary polyhedrons may profit even more, as the computational effort becomes even larger for general, especially non-convex polyhedra, cf.~\reffig{illustration_shaskov}. The performance gain achievable with the simple cuboid cells in FS3D is expected to provide a lower limit of the possible gain. %

\vspace{\baselineskip}

The following sections \ref{sec:strategy} and \ref{sec:computational_details} present the method in a general form, applicable to any polyhedron cell. %
A performance comparison with an accelerated bisection and decomposition-based approach is shown in a stand-alone PLIC positioning implementation as well as embedded into the flow solver FS3D. %
Originally, the flow solver FS3D employed the mentioned accelerated bisection approach, cf.\ appendix \ref{app:decomposition_approach}, which is used as a performance reference to compare the results of the enhanced three-phase PLIC positioning developed in the present study. %
The design of the numerical experiments and results are presented in section~\ref{sec:numerical_experiments}. %
A stand-alone evaluation and a prototypical three-phase test case simulated in FS3D show the potential of the new sequential three-phase positioning method both alone and within a parallelized flow solver. %


%

%
%
\subsection{Problem formulation}\label{subsec:problem_formulation}%
\begin{figure}[!htb]
\centering
\def\svgwidth{7cm}
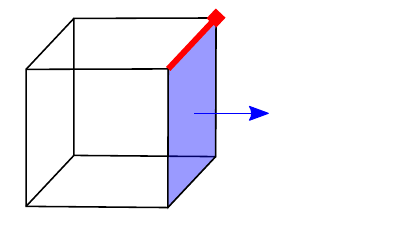
 \caption{Illustration of some important quantities of the polyhedron $\polyhedron*$, which is exemplary a cuboid here, but can be any polyhedron not intersecting itself. } \label{fig:Notation}
\end{figure}
Before we precisely state the formulation of the problem at hand, some notation associated with an arbitrary polyhedron $\polyhedron*\subset\setR^3$ representing a mesh cell is introduced. %
 The boundary $\partial\polyhedron*=\bigcup{\vfface_k}$ of the polyhedron $\polyhedron*$ is composed of $N^\vfface$ planar\footnote{Planarity is required for the application of the \textsc{Gaussian} divergence theorem, i.e., non-planar polygonal faces must be decomposed accordingly.} polygonal faces $\vfface_k$, each with an outer normal $\vn^\vfface_k$. %
Throughout this work, the subscript $k$ ($m$) refers to the faces (edges/vertices) of the respective polyhedron. %
Although the faces are allowed to be non-convex, which might occur for more than three edges on a face, we assume that the respective boundary $\partial\vfface_k$ admits no self-intersection. Note that connectivity between the faces is not required. %
Each face consists of the vertices $\vx^{\vfface}_{k,m}$ ($m=1\dots N^{\vfvert}_k$), which are ordered counter-clockwise with respect to $\vn^\vfface_k$. %
The respective summation limits, i.e., the number of faces and edges/vertices on a face, are suppressed for better readability. Also for notational convenience, the indices $m$ are assumed periodic, i.e., $\vx^{\vfface}_{k,N^{\vfvert}_k+1}=\vx^{\vfface}_{k,1}$. %
\refFig{Notation} illustrates the mentioned quantities. %
\begin{figure*}[!tb]
\null\hfill%
\includegraphics[page=3]{\thisfigurepath face_bracket_illustration}%
\hfill%
\includegraphics[page=1]{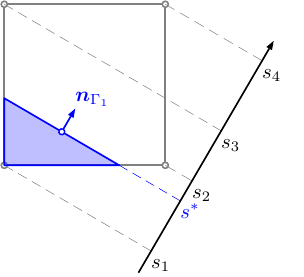}%
\hfill%
\includegraphics[page=2]{\thisfigurepath face_bracket_illustration}%
\hfill\null%
\caption{A prototypical three-phase PLIC configuration on a face (left):$\plicintersectionpoints$ denotes the intersection of the PLIC planes. This intersection is a triple line, thus a point on the face, if the intersection is within the cell and intersects the respective face. The original (center) and truncated (right) polygonal face are shown with signed distances $\primarysigndist_i$ and $\secondarysigndist_i$ of the vertices and $\primarysigndistref$ and $\secondarysigndistref$ of the resulting PLIC planes.}%
\label{fig:vertex_normal_projections}%
\end{figure*}
Next, we introduce the notation associated to the PLIC planes: The primary PLIC plane, employed to approximate the two interfaces of the primary phase towards the secondary phase and the ambient phase admits the unit normal $\primaryplicnormal$, as the sequential PLIC assumes no kink in the two interfaces. The secondary PLIC plane, which approximates the interface between the secondary phase and the ambient phase, has the unit normal $\secondaryplicnormal$. %
The unit normal orientations are exemplary shown in \reffig{vertex_normal_projections}.
With an arbitrary but fixed (cell-local) origin $\xbase$, the associated signed distances $\primarysigndist$ (primary) and $\secondarysigndist$ (secondary), respectively, define a family of PLIC planes as %
\begin{equation}
\begin{aligned}
& \primaryplicplane\fof{\primarysigndist}\!=\!\set{\vx\in\setR^3:\pliclvlset{1}\fof{\vx}=\primarysigndist}  %
\quad\text{and} \quad%
\\
& \secondaryplicplane\fof{\secondarysigndist}\!=\!\set{\vx\in\setR^3:\pliclvlset{2}\fof{\vx}=\secondarysigndist} %
\\
\quad\text{with}\quad & %
\pliclvlset{i}\fof{\vx}\!=\!\iprod{\vx-\xbase}{\vn_{\plicplane_i}} = (\vx-\xbase) \cdot \vn_{\plicplane_i}. 
\end{aligned}
\label{eqn:plic_parametrization_levelset}%
\end{equation}
The arbitrary (cell-local) origin $\xbase$ can be understood as a translation of the coordinate system into a cell-local one, the origin of which is typically chosen as a corner of the polyhedron cell. This way it is ensured that the computed distances within the polyhedron do not result in a significant precision loss due to large coordinate values. The primary $\primarysigndist$ and secondary $\secondarysigndist$ signed distances refer to the local origin $\xbase$. %
The associated negative and positive half-spaces of the family of PLIC planes are%
\begin{equation*}
\begin{aligned}
\neghalfspace[i]{u}\defeq\set{\vx\in\setR^3:\pliclvlset{i}\fof{\vx}\leq u} &%
\quad\text{and}\quad%
\\
\poshalfspace[i]{u}\defeq\set{\vx\in\setR^3:\pliclvlset{i}\fof{\vx}\geq u}. &%
\end{aligned}
\end{equation*}
\begin{figure*}[!tb]%
\null\hfill%
\includegraphics[page=1]{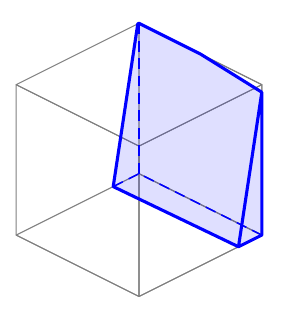}%
\hfill%
\includegraphics[page=2]{\thisfigurepath sequential_illustration}%
\hfill%
\includegraphics[page=3]{\thisfigurepath sequential_illustration}%
\hfill\null%
\caption{Illustration of sequential truncation for a cuboid: after the positioning of the primary plane $\primaryplicplane$, which truncates the primary volume (blue) from the polyhedron, one obtains the truncated residual polyhedron $\cutpolyhedron$ (center) The secondary plane $\secondaryplicplane$ (red) is positioned in the residual polyhedron $\cutpolyhedron$ truncating the secondary volume from the cut polyhedron. In this specific instance, the intersection of the PLIC planes, here a triple line denoted as $\plicintersection$ (magenta) intersects the faces $\vfface_4$ (front right face) and $\vfface_5$ (back left face) of the original polyhedron.}%
\label{fig:sequential_truncation}%
\end{figure*}%
We can now precisely state the problem to be solved: Given two pairs $\set{(\plicnormal,\refvof)}_i$ of unit normals $\vn_{\plicplane_i}$ and volume fractions $0<\primaryrefvof<1$ and $0<\secondaryrefvof<1-\primaryrefvof$ (both of which refer to the \textit{original} volume $\abs{\polyhedron*}$ of the polyhedron), signed distances $\primarysigndistref$ and $\secondarysigndistref$ of two sequentially positioned planes are sought in such a way that they truncate the given volume fractions from the polyhedral cell. This should be performed sequentially, in a primary positioning step in the original polyhedron $\polyhedron*$ and afterwards in a secondary positioning step in the residual polyhedron $\cutpolyhedron$. More precisely, the two positioning steps solve the following sub-problems:
\\
\textbf{primary positioning:} Find the \textit{primary} signed distance $\primarysigndistref$ such that $\abs{\polyhedron*\cap\neghalfspace[1]{\primarysigndistref}}=\primaryrefvof\abs{\polyhedron*}$, i.e., such that the plane $\primaryplicplane\fof{\primarysigndistref}$ truncates a sub-polyhedron of the volume $\primaryrefvof\abs{\polyhedron*}$. %
Let $\cutpolyhedron\defeq\polyhedron*\cap\poshalfspace[1]{\primarysigndistref}$ denote the residual sub-polyhedron with the volume $\brackets{1-\primaryrefvof}\abs{\polyhedron*}$, within which the secondary plane is positioned. %
\\
\textbf{secondary positioning:} Find the \textit{secondary} signed distance $\secondarysigndistref$ such that $\abs{\cutpolyhedron\cap\neghalfspace[2]{\secondarysigndistref}}=\secondaryrefvof\abs{\polyhedron*}$, i.e., such that the secondary plane $\secondaryplicplane$ truncates a sub-polyhedron with volume $\secondaryrefvof\abs{\polyhedron*}=\frac{\secondaryrefvof}{1-\primaryrefvof}\abs{\cutpolyhedron}$ from the residual polyhedron $\cutpolyhedron$. %
\\
\refFig{sequential_truncation} illustrates the sequential truncation procedure, showing the intersection line $\plicintersectionpoints$ of the two planes and the points $\vx_{0,k}$ of this line and the respective face of a prototypical cuboidal polyhedron.  %
Assuming, for the moment, that the PLIC planes are not parallel (i.e.,$\abs{\iprod{\primaryplicnormal}{\secondaryplicnormal}}\neq1$),  the line of intersection of the planes exists, and we denote $\plicintersectionpoints\defeq\primaryplicplane\fof{\primarysigndist}\cap\secondaryplicplane\fof{\secondarysigndist}$ as this intersection. The intersection is shown in three dimensions in \reffig{sequential_truncation} (right) and on a single face in \reffig{vertex_normal_projections} (left). Physically, this intersection represents an approximation of the so-called triple or contact line, if the intersection lies inside the cell. If the configuration is non- or fully wetted, this line lies outside the polyhedral cell. The positions of the sought PLIC planes are given by the signed distances $\primarysigndistref$ and $\secondarysigndistref$.

Let $\vec{X}$ and $\vec{X}^{\mathrm{cut}}$ be the set of vertices of the original ($\polyhedron*$) and truncated ($\cutpolyhedron$) polyhedron, respectively. The associated signed distances %
\begin{equation}
\begin{aligned}
\hat{\mathcal{S}}=\set*{\iprod{\vx-\xbase}{\primaryplicnormal}:\vx\in\vec{X}} &%
\quad\text{and}\quad%
\\
\hat{\mathcal{S}}^{\mathrm{cut}}=\set*{\iprod{\vx-\xbase}{\secondaryplicnormal}:\vx\in\vec{X}^{\mathrm{cut}}} &%
\end{aligned}
\label{eqn:vertex_normal_projections}%
\end{equation}
correspond to the positions $\primarysigndist$ and $\secondarysigndist$ for which the planes $\primaryplicplane$ and $\secondaryplicplane$ pass through a vertex. %
Note that it is possible that the plane passes through multiple vertices simultaneously. %
Let the elements $\projvertpos{i}$ of $\hat{\mathcal{S}}$ be arranged in ascending order ($\projvertpos{i}<\projvertpos{i+1}$) and define mutually disjoint so-called \textit{brackets} $\mathcal{B}_{i}\defeq(\projvertpos{i},\projvertpos{i+1}]$, %
starting with $\mathcal{B}_1\defeq[s_1,s_2]$, where the notation analogously applies to $\hat{\mathcal{S}}^{\mathrm{cut}}$, i.e., to the signed distances of the corners $\secondarysigndist_i$ of the residual polyhedron; %
cf.\ \reffig{vertex_normal_projections} for an illustration in two spatial dimensions. %
Within each bracket, the volume of a truncated polyhedron can be parametrized by a cubic polynomial in $\primarysigndist$ and $\secondarysigndist$, respectively. Note that, globally (i.e., on $\bigcup_i\mathcal{B}_i$ and $\bigcup_i\mathcal{B}^{\mathrm{cut}}_i$), the volume is a strictly monotone and, at least, continuous\footnote{An in-depth analysis of the regularity can be found in \cite{JCP_2021_fbip} in note 2.1 there.} function of the parameter$\primarysigndist$ or $\secondarysigndist$. %
The \textit{target} brackets $\refbracket$ and $\cutrefbracket$, respectively, contain the sought roots $\primarysigndistref$ and $\secondarysigndistref$, corresponding to the positions of the resulting PLIC planes. %

The positioning problem comprises two conceptually distinct subtasks: the efficient \textbf{computation of the volume} of truncated polyhedra as well as the \textbf{root-finding} of the scalar volume function. %
The present paper resorts to the root-finding algorithm of \citet{JCP_2021_fbip}, which was shown to be very efficient for several classes of polyhedra and an extensive set of combinations of volume fractions and normals; cf.\ \refapp{implicit_bracketing} for a description. %
The necessary adaptations of the volume computation for the sequential three-phase PLIC problem are described below: section~\ref{sec:strategy} outlines the general mathematical strategy for sequential truncation, while section~\ref{sec:computational_details} provides the details for the numerical treatment. Before, however, it is instructive to consider the possible topological phase configurations of the input data $\set{(\refvof,\plicnormal)}_i$ and their physical interpretation. %
\paragraph{Topological configurations} The sequential three-phase positioning may encounter different configurations as illustrated in \reffig{degenerate_configurations}, which are computed from the input data at runtime: %
\begin{figure*}[!tb]%
\null\hfill%
\includegraphics[page=4]{\thisfigurepath face_configuration}%
\hfill%
\includegraphics[page=5]{\thisfigurepath face_configuration}%
\hfill%
\includegraphics[page=2]{\thisfigurepath face_configuration}%
\hfill%
\includegraphics[page=3]{\thisfigurepath face_configuration}%
\hfill%
\includegraphics[page=1]{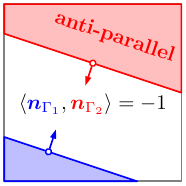}%
\hfill\null%
\caption{Topological classification of three-phase intersection configurations: triple, fully wetted and non-wetted, with associated degenerate cases for $\protect\iprod{\primaryplicnormal}{\secondaryplicnormal}=\pm1$. (For illustration purposes, we consider two spatial dimensions. However, the three-dimensional analoga are readily derived.)}%
\label{fig:degenerate_configurations}%
\end{figure*}%
\\
\textbf{triple:} In a triple configuration, the intersection $\plicintersectionpoints$ of the planes is partially contained within the convex hull of the cell, i.e., $\plicintersectionpoints\cap\convhull{\polyhedron*}\neq\emptyset$. Note that, for non-convex polyhedra, this does not necessarily correspond to a ``physical'' contact line, since $\plicintersectionpoints\cap\polyhedron*$ may be empty; cf.\ the right panel in \reffig{face_brackets}. %
\\
\textbf{fully wetted:} As the name suggests, a configuration of this type features a primary plane whose intersection with the polyhedron is fully contained in the negative half-space of the secondary plane, i.e., $\brackets{\primaryplicplane\cap\polyhedron*}\subset\neghalfspace[2]{}$. %
\\
\textbf{fully wetted (parallel degenerate, $\iprod{\primaryplicnormal}{\secondaryplicnormal}=1$):} The planes $\primaryplicplane$ and $\secondaryplicplane$ are aligned, implying that the negative half-space of $\primaryplicplane$ is a subset of the negative half-space of $\secondaryplicplane$. Thus, the secondary positioning problem corresponds to a translated version of the first positioning problem, i.e., $\set{\secondaryrefvof,\secondaryplicnormal,\cutpolyhedron}\equiv\set{\primaryrefvof+\secondaryrefvof,\secondaryplicnormal,\polyhedron*}$. %
\\
\textbf{non-wetted:} If the negative half-spaces of $\primaryplicplane$ and $\secondaryplicplane$ admit no overlap with the polyhedron, i.e., if $\brackets{\neghalfspace[1]{}\cap\neghalfspace[2]{}}\cap\polyhedron*=\emptyset$, the configuration is denoted non-wetted. %
\\
\textbf{non-wetted (anti-parallel degenerate, $\iprod{\primaryplicnormal}{\secondaryplicnormal}=-1$):} The negative half-spaces of $\primaryplicplane$ and $\secondaryplicplane$ do no overlap at all, implying that the primary truncation of the polyhedron $\polyhedron*$ does not affect the secondary truncation. %
Thus, $\set{\secondaryrefvof,\secondaryplicnormal,\cutpolyhedron}\equiv\set{\secondaryrefvof,\secondaryplicnormal,\polyhedron*}$, and one may solve $\set{\primaryrefvof,\primaryplicnormal,\polyhedron*}$ and $\set{\secondaryrefvof,\secondaryplicnormal,\polyhedron*}$ independently.

\begin{remark}[Ambiguity of degeneration]%
It is worth noting that the above defined topological configurations are not always unique. As can be seen in the center panel in \reffig{degenerate_configurations}, a parallel orientation, i.e., $\iprod{\primaryplicnormal}{\secondaryplicnormal}=1$, is not a necessary condition for $\set{\secondaryrefvof,\secondaryplicnormal,\cutpolyhedron}\equiv\set{\primaryrefvof+\secondaryrefvof,\secondaryplicnormal,\polyhedron*}$. %
In fact, one requires $\brackets{\primaryplicplane\cap\polyhedron*}\subset\brackets{\cutpolyhedron\cap\neghalfspace[2]{}}$ for a \textit{fully wetted} configuration, which may correspond to any  $\iprod{\primaryplicnormal}{\secondaryplicnormal}\in[-1,1]$. %

While such a more general definition obviously reduces the computational effort, its robust detection based on the input data becomes considerably involved, especially for non-convex polyhedra.
\end{remark}%

Identifying the aforementioned (anti-)parallel degenerated configurations from the input data in a first step reduces the number of exceptions within the algorithm at a later stage. Therefore, the remainder of this paper focusses on the non-degenerate cases, i.e., $\abs{\iprod{\primaryplicnormal}{\secondaryplicnormal}}<1$. In other words, in what follows, the primary and secondary plane are assumed to admit a non-empty intersection $\plicintersectionpoints=\primaryplicplane\cap\secondaryplicplane$, lying in- or outside the polyhedral cell. %
The (anti-)parallel configurations allow an independent positioning in three-phase cells without connectivity computations of the cut polyhedron anyway. %
The intersection of the primary and secondary truncating plane is the basis for the enhanced efficient volume computation of a twice truncated polyhedron.
%
%
%
%
%

%
%
\section{Volume computation for a twice truncated polyhedron}\label{sec:strategy}%
The efficient computation of the volume of the polyhedron resulting from the truncation of the original polyhedral cell with a single halfspace, employing the method introduced in \citet{JCP_2021_fbip}, allows for a major conceptual improvement when applied to the sequential truncation with two planes and the resulting volumes' computation in a three-phase cell: an appropriate choice of a reference point for the volume computation allows to skip the costly establishment of the cut polyhedron's connectivity, i.e.\ establishing the neighbor relations of points on the new face resulting from the truncation becomes unnecessary. %
Advantageously, the analytic considerations concerning the regularity derived for the two-phase case directly apply to the sequential truncation. Hence, instead of reproducing them here, the interested reader is referred to \cite[subsection~2.1]{JCP_2021_fbip}.
By application of the \textsc{Gaussian} divergence theorem, the volume of a truncated polyhedron $\polyhedron*$ with outer unit normal $\vn_{\polyhedron*}\in\set{\vn^\vfface_k}$ can be expressed as %
\begin{equation}
\begin{aligned}
\abs{\neghalfspace[1]{}\cap\polyhedron*}
& =%
\int\limits_{\neghalfspace[1]{}\cap\polyhedron*}{1\dvol}
=%
\int\limits_{\neghalfspace[1]{}\cap\polyhedron*}{\frac{1}{3}\nabla \cdot \brackets{\vx-\plicbase}\dvol} 
\\
& =%
\frac{1}{3}\int\limits_{\partial\brackets{\neghalfspace[1]{}\cap\polyhedron*}}{\iprod{\vx-\plicbase}{\vn_{\polyhedron*}}\darea(\vx)}
\\%
& = 
\frac{1}{3}\Biggl(\sum\limits_{k}{\iprod{\vx^\vfface_{k,1}-\plicbase}{\vn^\vfface_k}\abs{\neghalfspace[1]{}\cap\vfface_k}}%
\\
& + \underbrace{\brackets*{\iprod{\xbase-\plicbase}{\primaryplicnormal}+\primarysigndist}}_{= 0 \text{ for }\plicbase \text{ coplanar to } \primaryplicplane} \abs{\primaryplicplane\cap\polyhedron*}\Biggl),%
\end{aligned}
\label{eqn:volume_computation_single}
\end{equation}
where $\dd{o}$ denotes the two-dimensional \textsc{Lebesgue} measure (area) and $\plicbase$ is an arbitrary but constant point. Note that some constraints on $\plicbase$ appear in the following, which are required to exploit the advantages of the sequential truncation's shared geometry. \refEqn{volume_computation_single} casts the immersed volume as a sum of face-based quantities, each of which is a binary product of the immersed area and the signed distance to the reference point $\plicbase$, which is fixed by some appropriate choices in the following: The first constraint on $\plicbase$ arises from the fact that one of the boundary segments is located on the truncating plane $\primaryplicplane\cap\polyhedron*$. %
Its contribution to \refeqn{volume_computation_single} becomes zero, if one chooses any reference point $\plicbase$ coplanar~\footnote{This means that $\set{\plicbase}\cup\primaryplicplane$ forms a coplanar set of points, while $\plicbase\in\primaryplicplane\cap\polyhedron*$ is not necessarily required. The reference point $\plicbase$ can be outside the polyhedron, e.g., if a non-wetted or fully-wetted configuration is present.} to $\primaryplicplane$, i.e., if $\iprod{\xbase-\plicbase}{\primaryplicnormal}+\primarysigndist=0$. 

Analogously to \refeqn{volume_computation_single}, the volume of a sequentially truncated residual polyhedron reads %
\begin{equation}
\begin{aligned}
\abs{\neghalfspace[2]{}\cap\cutpolyhedron}%
=&  \frac{1}{3}\Biggl(\sum\limits_{k}{\iprod{\vx^\vfface_k-\plicbase}{\vn^\vfface_k}\abs{\neghalfspace[2]{}\cap\poshalfspace[1]{}\cap\vfface_k}}%
\\
  & + \underbracea{\underbrace{\brackets*{\iprod{\xbase-\plicbase}{\primaryplicnormal}+\primarysigndist}}_{= 0 \text{ for }\plicbase \text{ coplanar to } \primaryplicplane}\abs{\primaryplicplane\cap\neghalfspace[2]{}\cap\polyhedron*}}%
\\
 & \underbracebd{+ \underbrace{\brackets*{\iprod{\xbase-\plicbase}{\secondaryplicnormal}+\secondarysigndist}}_{= 0 \text{ for }\plicbase \text{ coplanar to } \secondaryplicplane}\abs{\secondaryplicplane\cap\cutpolyhedron}\vphantom{\sum\limits_k}}_{= 0 \text{ for }\plicbase \in \plicintersectionpoints= \primaryplicplane\cap\secondaryplicplane}\Biggr),%
\end{aligned}
\label{eqn:volume_computation_sequential}%
\end{equation}
where the boundary contains two segments on the truncating planes, $\primaryplicplane\cap\neghalfspace[2]{}\cap\polyhedron*$ on the primary and $\secondaryplicplane\cap\cutpolyhedron$ on the secondary truncating plane. %
Choosing the reference point $\plicbase$ introduced in \refeqn{volume_computation_single} coplanar to the primary PLIC plane $\primaryplicplane$ eliminates the second term on the right-hand side of \refeqn{volume_computation_sequential}. In order to also eliminate the last term in \refeqn{volume_computation_sequential}, the reference point $\plicbase$ needs to be coplanar to both truncating planes, $\primaryplicplane\fof{\primarysigndist}$ and $\secondaryplicplane\fof{\secondarysigndist}$, corresponding to
\begin{equation}
\begin{aligned}
 \iprod{\xbase-\plicbase}{\primaryplicnormal}+\primarysigndist & =0\text{ and }\\ 
 \iprod{\xbase-\plicbase}{\secondaryplicnormal}+\secondarysigndist & =0 \text{.} \\ %
\end{aligned} \label{eqn:volumeoriginconstraints}
\end{equation}

Thus, we further restrict the choice of $\plicbase$ to lie on the intersection $\plicintersectionpoints = \primaryplicplane\cap\secondaryplicplane$ of the two planes; cf. \reffig{sequential_truncation}~(right) for an illustration.%

Let $\cutvfface_{k}\defeq\poshalfspace[1]{}\cap\vfface_k$ denote the truncated version of the original face $\vfface_k$. %
The immersed face area of the first truncation $\abs{\neghalfspace[1]{}\cap\vfface_{k}}$ is computed analogous to \citet{JCP_2021_fbip}, but with the same choice of reference points described in the following, which is advantageous for the computation of the secondary immersed face area $\abs{\neghalfspace[2]{}\cap\cutvfface_{k}}$. %
We assign the outer co-normal $\vN_{k,m}$ with $\iprod{\vn^\vfface_k}{\vN_{k,m}}~=~0$ to the edges $\vfedge_{k,m}=\convhull{\vx^{\vfface}_{k,m},\vx^{\vfface}_{k,m+1}}$.  $\convhull{...}$ is defined as the convex hull of the given points. %
With the vertices sorted counter-clockwise, $\vN_{k,m}$ can be computed by normalizing $\cross{\brackets{\vx^{\vfface}_{k,m+1}-\vx^{\vfface}_{k,m}}}{\vn^\vfface_k}$ to become a unit vector. %
By recursive application of the \textsc{Gaussian} divergence theorem, in a manner analogous to the volume computation described above, the secondary immersed area $\abs{\neghalfspace[2]{}\cap\cutvfface_{k}}$ on the original face $\vfface_k$ becomes %
\colorlet{pcol}{black}
\colorlet{scol}{black}
\begin{equation}
\begin{aligned}
\abs{\neghalfspace[2]{}\cap & \cutvfface_{k}}%
\\
 = & \frac{1}{2}\sum\limits_{m}{\iprod{\vx^\vfface_{k,m}-\vx_{0,k}}{\vN_{k,m}}\abs{\neghalfspace[2]{}\cap\poshalfspace[1]{}\cap\vfedge_{k,m}}} & \\%
 + & \frac{1}{2}\underbrace{\iprod*{\vy_{\primaryplicplane,k}-\vx_{0,k}}{\frac{\primaryplicnormal-\iprod{\primaryplicnormal}{\vn^\vfface_k}\vn^\vfface_k}{1-\iprod{\primaryplicnormal}{\vn^\vfface_k}^2}}}_{=0 \text{ for } \vx_{0,k} \text{ coplanar to } \primaryplicplane \text{ and } \vfface_k}%
\abs{\primaryplicplane\cap\neghalfspace[2]{}\cap\vfface_k}\\%
 + & \frac{1}{2}\underbrace{\iprod*{\vy_{\secondaryplicplane,k}-\vx_{0,k}}{\frac{\secondaryplicnormal-\iprod{\secondaryplicnormal}{\vn^\vfface_k}\vn^\vfface_k}{1-\iprod{\secondaryplicnormal}{\vn^\vfface_k}^2}}}_{=0 \text{ for } \vx_{0,k} \text{ coplanar to } \secondaryplicplane \text{ and }  \vfface_k}%
\abs{\secondaryplicplane\cap\cutvfface_{k}}%
\end{aligned}
\label{eqn:face_area_sequential}%
\end{equation}
with arbitrary points $\vy_{\plicplane_i,k}$ coplanar to both $\vfface_k$ and $\plicplane_i$ (\textcolor{blue}{$\blacksquare$}/\textcolor{red}{$\blacksquare$} in \reffig{faceorigin_illustration} illustrate some candidates). %
In order to eliminate the last two contributions in \refeqn{face_area_sequential} one finds, after some manipulations, that the reference point for the area computation $\vx_{0,k}$ associated to the face $\vfface_k$ is subject to three scalar constraints: %
\begin{figure*}[!tb]%
\null\hfill%
\includegraphics[page=2]{\thisfigurepath origin_illustration}%
\hfill%
\includegraphics[page=1]{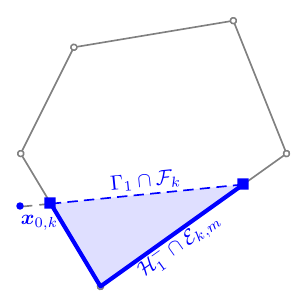}%
\hfill%
\includegraphics[page=3]{\thisfigurepath origin_illustration}
\hfill\null%
\caption{Computation of the area of face $\vfface_k$ immersed in the negative half-space of the primary (blue) and secondary (red) plane. %
The bold lines (center/right) indicate the boundary segments with non-zero contribution to eqs.~(\ref{eqn:volume_computation_single}) and (\ref{eqn:volume_computation_sequential}), respectively, %
assuming the depicted reference points on the polyhedrons faces $\vx_{0,k}$ (\textcolor{blue}{$\bullet$}/\textcolor{magenta}{$\bullet$}). %
It becomes clear that for each truncated face $\vfface_k$ only one choice of $\vx_{0,k}$ remains which allows the area computation (and thus also the volume computation) without reestablishing connectivity: %
The last constraint in eq.~(\ref{eqn:face_origin_constraints}) drops out for the primary truncation (center), for any choice of $\vx_{0,k}$ on the blue dashed line. %
In the presence of a secondary plane (right), the choice of $\vx_{0,k}$ which simultaneously eliminates both of the planes' contributions is unique on the intersection of the red and blue dashed lines.}%
\label{fig:faceorigin_illustration}%
\end{figure*}%
\begin{equation}
\begin{aligned}
\iprod{\vx_{0,k}-\vx^\vfface_{k,1}}{\vn^\vfface_k}=0&\quad\text{(coplanar to $\vfface_k$)},\\%
\iprod{\vx_{0,k}-\xbase}{\primaryplicnormal}=\primarysigndist&\quad\text{(coplanar to $\primaryplicplane\fof{\primarysigndist}$)},\\%
\iprod{\vx_{0,k}-\xbase}{\secondaryplicnormal}=\secondarysigndist&\quad\text{(coplanar to $\secondaryplicplane\fof{\secondarysigndist}$)},%
\end{aligned}
\label{eqn:face_origin_constraints}%
\end{equation}
corresponding to a \textbf{unique} choice of $\vx_{0,k}$ on the intersection $\plicintersectionpoints$ and on the respective face $\vfface_k$; cf.~the right panel in \reffig{faceorigin_illustration} for an illustration. Note that this choice of the reference point $\vx_{0,k}$ ensures that the contributions of the edges located on the planes $\primaryplicplane$ and $\secondaryplicplane$ emerging from the truncation become zero. The choice of the volume reference point $\plicbase$ on the intersection $\plicintersectionpoints = \primaryplicplane\cap\secondaryplicplane$ and the area reference points $\vx_{0,k}$ on $\plicintersectionpoints$ and additionally on face $\vfface_k$ eliminates the necessity of reestablishing the connectivity of the residual polyhedron $\cutpolyhedron$ after the first truncation, i.e.\ the neighbor relations of points and faces on the new face of the cut polyhedron after the truncation are not required. %
Since any area reference point $\vx_0,k$ is an admissible candidate for the volume reference point $\plicbase$, we set $\plicbase\defeq \vx_{0,1}$ (on the first \textit{intersected} face chosen as $\vfface_1$), which reduces the computation effort further by additionally eliminating the contribution of $\vfface_1$ to the volume. %
\begin{remark}[Multiple truncation]%
A secondary implication of \refeqn{face_origin_constraints} is that the \textsc{Gaussian} divergence theorem cannot be employed to generalize this concept for positioning more than two planes sequentially, since the choice of $\vx_{0,k}\in\setR^3$ would be subject to more than three constraints. The topological connectivity of the truncated polyhedron must be reestablished $\lfloor\frac{M-1}{2}\rfloor$ times for a multi-material application with $M$ different materials.%
\end{remark}%

In the upcoming section~\ref{sec:computational_details} we derive the computational details of the strategy outlined above. %
%

%
%
\section{Computational details}\label{sec:computational_details}%
\begin{figure*}[!tb]%
\null\hfill%
\includegraphics[page=1]{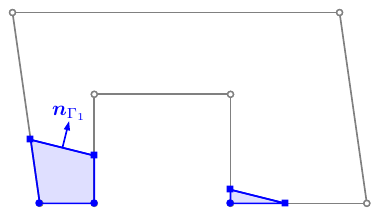}%
\hfill%
\includegraphics[page=2]{\thisfigurepath face_cutface}%
\hfill\null%
\caption{Example of a non-convex cell face: Original face $\vfface_k$ (left) and associated truncation $\cutvfface_k$ (right), defined as a list of truncated edges $\cutvfedge_{k,m}$. The $\textcolor{gray}{\circ}$ and $\textcolor{blue}{\blacksquare}$ indicate the vertices of the original polygon and the intersections with $\primaryplicplane$, respectively. Note that the edges $\cutvfedge_{k,m}$ do not form a closed polygon, because the boundary segments located on the primary plane $\primaryplicplane$ (dashed) are not contained. Due to the choice of the origin $\vx_{0,k}$, their contributions to the area sum become zero anyway.}%
\label{fig:face_truncated_edges}%
\end{figure*}

%
%
%
\begin{table}[!bt]
\centering%
\caption{Edge status as function of the status of the associated vertices $\vx^\vfface_{k,m}\protect\defeq*\vfvert_u$ and $\vx^\vfface_{k,m+1}\protect\defeq*\vfvert_v$; cf.~\refnote{floating_point_operations}.}%
\label{tab:intersection_status_edge}%
\begin{tabular}{c||ccc|cc|cc|c}
&\multicolumn{3}{c|}{}&\multicolumn{5}{c}{\textbf{degenerate}}\\%
&\rotatebox{90}{exterior}&\rotatebox{90}{interior}&\rotatebox{90}{intersected}&\multicolumn{2}{c|}{\rotatebox{90}{exterior}}&\multicolumn{2}{c|}{\rotatebox{90}{interior}}&\rotatebox{90}{intersected}\\%
&\includegraphics[page=4]{\thisfigurepath edgestatus_illustration}&\includegraphics[page=2]{\thisfigurepath edgestatus_illustration}&\includegraphics[page=1]{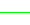}&\multicolumn{2}{c|}{\includegraphics[page=5]{\thisfigurepath edgestatus_illustration}}&\multicolumn{2}{c|}{\includegraphics[page=3]{\thisfigurepath edgestatus_illustration}}&\includegraphics[page=6]{\thisfigurepath edgestatus_illustration}\\%
\hline%
$\vfstatus[i]{\vfvert_u}$&$1$&$-1$&$\pm1$&$1$&$0$&$-1$&$\phantom{-}0$&$0$\\%
$\vfstatus[i]{\vfvert_v}$&$1$&$-1$&$\mp1$&$0$&$1$&$\phantom{-}0$&$-1$&$0$\\%
\hline%
$\vfstatus[i]{\vfedge}$&$1$&$-1$&$0$&\multicolumn{2}{c|}{$2$}&\multicolumn{2}{c|}{$-2$}&$3$%
\end{tabular}%
\\
\begin{tabular}{c}
\includegraphics[page=7]{\thisfigurepath edgestatus_illustration}%
\end{tabular}
\end{table}%
The robust treatment of intersections requires a hierarchically consistent evaluation of the topological properties of vertices $\set{\vx^\vfface_{k,m}}$ and edges $\set{\vfedge_{k,m}}$. With respect to the oriented hypersurface $\plicplane_i$, described by the level-set function given in \refeqn{plic_parametrization_levelset}, the logical status $\vfstatus[i]{}$ of any vertex is either interior ($\vfstatus[i]{}=-1$), intersected ($\vfstatus[i]{}=0$) or exterior ($\vfstatus[i]{}=1$). %
The hierarchically superior entity is an edge $\vfedge_{k,m}$, whose status is uniquely determined by the status of its associated vertices $\vx^\vfface_{k,m}$ and $\vx^\vfface_{k,m+1}$. %
An edge with at least one of its vertices located on the plane $\plicplane_i$ becomes degenerate, where the type of degeneration depends on the status of the second vertex of that edge; cf.~\reftab{intersection_status_edge}. %
\begin{note}[Floating point operations]\label{note:floating_point_operations}%
For the purpose of numerical robustness, the status assignment employs a tubular neighborhood of thickness $2\zerotol$ around $\plicplane_i$ corresponding to the interval $(-\zerotol,\zerotol)$. One obtains the logical status (with respect to the plane $\plicplane_i$) of a vertex $\vx$ as %
\begin{align*}%
\vfstatus[i]{\vx}\defeq%
\begin{cases}%
\phantom{-}0&\text{if}\quad\abs{\pliclvlset{i}\fof{\vx}}<\zerotolerance,\\%
\sign{\pliclvlset{i}\fof{\vx}}&\text{if}\quad\abs{\pliclvlset{i}\fof{\vx}}\geq\zerotolerance.%
\end{cases}%
\end{align*}
Thus, any point whose absolute distance to $\plicplane_i$ falls below $\zerotol$ is considered to be on $\plicplane_i$. Note that the choice of an appropriate tolerance strongly depends on the absolute value of the characteristic length of the polyhedron $\polyhedron*$. Throughout this work, we let $\zerotolerance\defeq\num{e-14}$. In fact, all zero-comparisons are implemented in this way. %
\end{note}%
As outlined in section~\ref{sec:strategy}, we employ the \textsc{Gaussian} divergence theorem to compute the volume of the intersection of a polyhedron $\polyhedron*$ and the negative half-space of a plane $\primaryplicplane$, cf. \refeqn{volume_computation_single}. %
In \refeqn{volume_computation_single} we applied the divergence theorem to convert the volume integral into a surface integral. %
With the choice of $\plicbase = \xbase + \primarysigndist \primaryplicnormal$ which is not only coplanar to the primary PLIC plane as required, but is additionally a convenient choice  as the vector is additionally orthogonal to the primary PLIC plane, one obtains the primary volume fraction dependent on the primary signed distance 
\begin{equation}
\begin{aligned}
\primarypolyvof\fof{\primarysigndist} & = \frac{\abs{\neghalfspace[1]{s}\cap\polyhedron*}}{\abs{\polyhedron*}}=\\
& \frac{1}{3\abs{\polyhedron*}}\sum\limits_{k}{%
\Biggl(\iprod{\vx^\vfface_{k,1}-\xbase}{\vn^\vfface_k} 
 -\primarysigndist\iprod{\primaryplicnormal}{\vn^\vfface_k}\Biggr)%
\primaryimmersedarea{k}\fof{\primarysigndist}}%
\text{,}%
\end{aligned}
\label{eqn:volume_fraction_primary}%
\end{equation}
where we write $\primaryimmersedarea{k}\defeq\abs{\neghalfspace[1]{}\cap\vfface_k}$; cf.~\refeqn{volume_computation_single}. %
Concerning \refeqn{volume_fraction_primary}, two points are worth noting: Firstly, applying the divergence theorem to transforms the volume integral into a surface integral, where, for a polyhedron, the integrand can be expressed as a sum of inner products involving the normals of the polyhedral faces. %
This removes any terms perpendicular to $\primaryplicnormal$. %
Secondly, in the first inner product $\vx^\vfface_{k,1}$ could be replaced by any point coplanar to the face $\vfface_k$. %
After finding the signed distance $\primarysigndistref$ corresponding to $\primaryrefvof$ with the method described in \citet{JCP_2021_fbip} and \refapp{implicit_bracketing}, one obtains a truncated residual polyhedron $\cutpolyhedron\fof{\primarysigndistref}$. %
We choose a representation in terms of truncated faces $\cutvfface_{k}$, each of which corresponds to a list of truncated edges $\cutvfedge_{k,m}\defeq\poshalfspace[1]{}\cap\vfedge_{k,m}$. %
It becomes clear in the following, why it is sufficient to only include the (truncated) edges of the original polyhedron. The contributions of the truncated faces are eliminated not only in the volume computation, but also for the area computation discussed below. %
This is a major advantage of our efficient positioning method as no corner to edge and edge to face connectivity relation has to be recomputed after the first truncation. %
The truncated edges $\cutvfedge_{k,m}$ are obtained based on the logical status (cf.~\reftab{intersection_status_edge}) of the vertices $\vx^\vfface_{k,m}$ and $\vx^\vfface_{k,m+1}$ of the original edge $\vfedge_{k,m}$: %
\begin{equation}%
\begin{aligned}%
\cutvfedge_{k,m}\defeq%
\begin{cases}%
\emptyset & \text{if the edge $\vfedge_{k,m}$ is } \\
         & \text{$\primaryplicplane$-interior, i.e., } \\
         & \vfstatus[1]{\vfedge_{k,m}}=-1,\\[10pt]%
\vfedge_{k,m}&\text{if the edge $\vfedge_{k,m}$ is} \\
         & \text{$\primaryplicplane$-exterior, i.e., } \\
         & \vfstatus[1]{\vfedge_{k,m}}=1,\\[10pt]%
\convhull{\vx^\vfface_{k,m},\hat{\vx}\fof{\primarysigndistref}}&\text{if the edge is intersected} \\
         & \text{ ($\vfstatus{\vfedge_{k,m}}=0$)} \\
         & \text{and\ }\vfstatus[1]{\vx^\vfface_{k,m}}=1,\\[10pt]%
\convhull{\hat{\vx}\fof{\primarysigndistref},\vx^\vfface_{k,m+1}}&\text{if the edge is intersected} \\
         & \text{($\vfstatus{\vfedge_{k,m}}=0$)} \\
         & \text{and\ }\vfstatus[1]{\vx^\vfface_{k,m}}=-1,\\%
\end{cases}%
\end{aligned}%
\label{eqn:truncated_edge_assignment}
\end{equation}%
where the intersection of the edge $\vfedge_{k,m}$ with the (primary) plane $\primaryplicplane$ contains precisely the point %
\begin{align}%
\hat{\vx}\fof{\primarysigndistref}=%
\vx^\vfface_{k,m}+%
\frac{\primarysigndistref+\iprod{\xbase-\vx^\vfface_{k,m}}{\primaryplicnormal}}{\iprod{\vx^\vfface_{k,m+1}%
-\vx^\vfface_{k,m}}{\primaryplicnormal}}\brackets*{\vx^\vfface_{k,m+1}-\vx^\vfface_{k,m}}.\label{eqn:truncated_edge_intersection}%
\end{align}%
\refEqn{truncated_edge_assignment} can be understood as follows: interior (exterior) edges are discarded (copied), while the $\primaryplicplane$-interior vertex ($\vfstatus[1]{}=-1$) of an intersected edge is replaced by the intersection with $\primaryplicplane$. \refFig{face_truncated_edges} provides an illustration. %
Note that the truncated edges do not need to form a closed polygon, since no edge connectivity information is required, due to the intended application of the \textsc{Gaussian} divergence theorem with the previously discussed advantageous choice of reference points $\plicbase\fof{\secondarysigndist;\primarysigndistref}$ for the volume- and $\vx_{0,k}\fof{\secondarysigndist;\primarysigndistref}$ for the area computations. %
The intersection $\plicintersectionpoints$ of the two planes is parametrized as %
\begin{equation}
\plicintersectionpar\fof{u;\primarysigndist,\secondarysigndist}=\vy_0\fof{s}+\secondarysigndist\vtau_\plicintersectionpar+u\vmu_\plicintersectionpar\label{eqn:intersection_line}%
\end{equation}
with %
\begin{equation}
\begin{aligned}
\vy_0\fof{s} & =  \xbase+\primarysigndist\frac{\primaryplicnormal-\normalip\secondaryplicnormal}{1-\normalip^2}  ,\quad \\%
\vtau_{\plicintersectionpar} & =  \frac{\secondaryplicnormal-\normalip\primaryplicnormal}{1-\normalip^2}  \quad\text{and}\quad \\
\vmu_{\plicintersectionpar} & =  \frac{\cross{\primaryplicnormal}{\secondaryplicnormal}}{1-\normalip^2}  .%
\end{aligned}
\label{eqn:intersection_line_definition}
\end{equation}
\refEqn{intersection_line} can be understood as follows: the signed distances of the PLIC planes, i.e., $\primarysigndist$ and $\secondarysigndist$, determine the intersection $\plicintersectionpoints$, the parametrization of which is $\plicintersectionpar$, running along the segment of the intersection of the two planes with $u$ as the curve parameter. %
In analogy to the truncation of the original polyhedron $\polyhedron*$, the secondary volume fraction is obtained by application of the \textsc{Gaussian} divergence theorem to the truncated polyhedron $\cutpolyhedron$ and choosing $\plicbase$ on $\plicintersection$, cf.\  \refeqn{intersection_line}, i.e., %
\begin{equation}
\begin{aligned}
\secondarypolyvof\fof{\secondarysigndist;\cutpolyhedron\fof{\primarysigndistref}}%
& =\frac{\abs{\neghalfspace[2]{\secondarysigndist}\cap\cutpolyhedron\fof{\primarysigndistref}}}{\abs{\polyhedron*}} \\%
& =\frac{1}{3\abs{\polyhedron*}}\sum\limits_{k}{%
\brackets{\volcoeffconst{k} - \secondarysigndist\volcoefflin{k}}%
\immersedarea{k}\fof{\secondarysigndist;\primarysigndistref}},%
\end{aligned}
\label{eqn:volume_fraction_secondary}%
\end{equation}
with the immersed area $\immersedarea{k}\defeq\abs{\neghalfspace[2]{}\cap\poshalfspace[1]{}\cap\vfface_k}$; cf.~\reffig{faceorigin_illustration}. %
The volume coefficients are set as %
\begin{align}
\volcoeffconst{k}\defeq\iprod{\vx^\vfface_{k,1}-\vy_0\fof{\primarysigndistref}}{\vn^\vfface_k}%
\quad\text{and}\quad%
\volcoefflin{k}\defeq\iprod{\vtau_\plicintersectionpar}{\vn^\vfface_k}\text{.}%
\label{eqn:area_coefficients}%
\end{align}
For reasons of computational robustness, the secondary volume fraction in \refeqn{volume_fraction_secondary} refers to the volume of the \textbf{truncated} polyhedron, i.e., $0\leq\secondarypolyvof\leq1$. Thus, the target secondary volume fraction, which refers to the original polyhedron cell, must be rescaled accordingly, i.e., $\secondaryrefvof$ becomes $\frac{\secondaryrefvof}{1-\primaryrefvof}$. %
To improve the readability of the notation, we drop the asterisk for the primary and secondary signed distance for the remainder of this paper, implying that $\primarysigndist$ from now on means $\primarysigndistref$ and $\secondarysigndist$ means $\secondarysigndistref$. %

%
%
\subsection{The computation of the immersed area $\immersedarea{k}$}%
\begin{figure*}[!tb]%
\null\hfill%
\includegraphics[page=1]{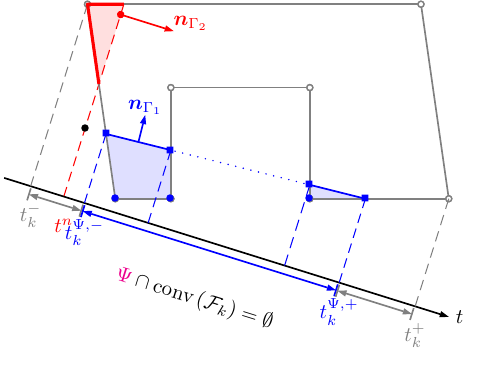}%
\hfill
\includegraphics[page=2]{\thisfigurepath face_illustration_quantities}%
\hfill\null%
\caption{Example of a non-convex cell face: Illustration of the limits of the signed distances on a polyhedrons face in eq.~(\ref{eqn:immersed_face_delimination}) for the computation of the immersed area $\immersedarea{k}$ of a truncated face $\cutvfface_k$ with different intersection configurations. %
In the left panel, the intersection of $\plicintersectionpoints$ with the plane containing $\vfface_k$ (\textcolor{black}{\textbullet}) is not located in its convex hull, implying that the original edges $\vfedge_{k,m}$ (bold) can be employed for the area computation (second case in \refeqnsafe{immersed_face_area}). %
In the right panel, with $\vx_{0,k}\in\protect\convhull{\vfface_k}$, the truncated edges $\cutvfedge_{k,m}$ (bold) are employed (first case in eq.~(\ref{eqn:immersed_face_area})).}%
\label{fig:face_brackets}%
\end{figure*}

As stated in \refeqn{face_area_sequential}, the computation of the immersed area exploits the \textsc{Gaussian} divergence theorem, i.e., %
\begin{align}
\immersedarea{k}=%
\frac{1}{2}\sum\limits_{m}{\iprod{\vx^\vfface_{k,m}-\vx_{0,k}}{\vN_{k,m}}\abs{\cutvfedge_{k,m}}}.\label{eqn:face_area_sequential_computation}%
\end{align}
For $\iprod{\vmu_\plicintersectionpar}{\vn^\vfface_k}\neq0$, i.e., if the triple line $\plicintersectionpoints$ is not parallel to the face, the reference point $\vx_{0,k}$ is the unique intersection of $\plicintersectionpoints$ with the plane containing the face $\vfface_k$, namely %
\begin{equation}
\begin{aligned}
\vx_{0,k}\fof{\secondarysigndist;\primarysigndist}= & \brackets*{\vy_0\fof{\primarysigndist}-\frac{\iprod{\vy_0\fof{\primarysigndist}%
-\vx^\vfface_{k,1}}{\vn^\vfface_k}}{\iprod{\vmu_\plicintersectionpar}{\vn^\vfface_k}}\vmu_\plicintersectionpar}
\\
& +%
\secondarysigndist\brackets*{\vtau_\plicintersectionpar-\frac{\iprod{\vtau_\plicintersectionpar}{\vn^\vfface_k}}{\iprod{\vmu_\plicintersectionpar}{\vn^\vfface_k}}\vmu_\plicintersectionpar} \\
\defeq* & \vy_{0,k}\fof{\primarysigndist}+\secondarysigndist \vtau_{\plicintersectionpar,k},%
\end{aligned}
\label{eqn:referencepointchoice}
\end{equation}
where \refeqn{intersection_line_definition} contains the definition of the involved quantities. %
In general, the intersection is not necessarily located in the convex hull of the face $\vfface_k$, i.e., $\vx_{0,k}\not\in\convhull{\vfface_k}$; cf.~\reffig{face_brackets}.
With the coefficients %
\begin{equation}
\begin{aligned}
\areacoeffconst{k}{m}=%
\begin{cases}%
\iprod{\vx^\vfface_{k,m}-\vy_{0,k}\fof{\primarysigndist}}{\vN_{k,m}}\text{ if }\secondarysigndist\in[\sdlowertriple{k},\sduppertriple{k}],\\[10pt]%
\iprod{\vx^\vfface_{k,m}}{\vN_{k,m}}+\frac{\iprod{\vx^\vfface_{k,1}}{\vn^\vfface_k}\iprod{\vn^\vfface_k}{\secondaryplicnormal}-\iprod{\xbase}{\secondaryplicnormal}}{1-\iprod{\vn^\vfface_k}{\secondaryplicnormal}^2}\iprod{\vN_{k,m}}{\secondaryplicnormal} \\
\quad \text{ if }\secondarysigndist\not\in[\sdlowertriple{k},\sduppertriple{k}],%
\end{cases}%
\end{aligned}%
\label{eqn:area_coefficients}
\end{equation}%

the secondary immersed area can be cast as %
\begin{equation}
\begin{aligned}
\immersedarea{k} = & %
\frac{1}{2}\sum\limits_{m}\Bigl(%
\bigl(\areacoeffconst{k}{m}\fof{\primarysigndist} \\
& +\secondarysigndist\areacoefflin{k}{m}\fof{\primarysigndist}\bigr)%
\brackets*[c]{\begin{matrix}
\ell\fof{\cutvfedge_{k,m};\secondarysigndist} \text{ if } \secondarysigndist\in[\sdlowertriple{k},\sduppertriple{k}]\\%
\ell\fof{\vfedge_{k,m};\secondarysigndist} \text{ if }\secondarysigndist\not\in[\sdlowertriple{k},\sduppertriple{k}]%
\end{matrix}}\Biggr) \\
& -%
\brackets*[c]{\begin{matrix}
\primaryimmersedarea{k} \text{ if }\secondarysigndist>\sduppertriple{k}\\%
0&\text{ if } \secondarysigndist\leq\sduppertriple{k}%
\end{matrix}}%
\end{aligned}
\label{eqn:immersed_face_area}%
\end{equation}
with the immersed edge length $\ell\fof{\vfedge;\secondarysigndist}\defeq\abs{\neghalfspace[2]{\secondarysigndist}\cap\poshalfspace[1]{\primarysigndistref}\cap\vfedge}$ and the parameters %
\begin{equation}
\begin{aligned}%
\sdlowerface{k} & \defeq\min_m\set{\pliclvlset{2}\fof{\vx^\vfface_{k,m}}:\pliclvlset{1}\fof{\vx^\vfface_{k,m}}\geq0}%
\quad\text{and}\quad \\%
\sdlowertriple{k} & \defeq\min\set{\pliclvlset{2}\fof{\vx}:\vx\in\primaryplicplane\cap\vfface_k},%
\end{aligned}%
\label{eqn:immersed_face_delimination}
\end{equation}
where $\sdupperface{k}$ and $\sduppertriple{k}$ are defined analogously with the respective maxima. See \reffig{face_brackets} for an illustration. %
The formulation of \refeqn{immersed_face_area} as a product allows to conveniently evaluate the derivatives of $\immersedarea{k}$ with respect to $\secondarysigndist$, which are required for the root-finding part of the positioning procedure; cf.~\refapp{implicit_bracketing}. %

\begin{remark}[Geometric interpretation]%
We denote the smallest (largest) value of the signed distance $\secondarysigndist$ for which the face $\vfface_k$ is intersected by the PLIC plane $\secondaryplicplane$ by $\sdlowerface{k}$ ($\sdupperface{k}$). %
For any $\secondarysigndist\in[\sdlowertriple{k},\sduppertriple{k}]$, the intersection $\vx_{0,k}\fof{\secondarysigndist}$ of the triple line $\plicintersectionpoints$ with the plane containing $\vfface_k$ is contained in the convex hull of $\vfface_k$. %
The case selection in \refeqs{area_coefficients} and \refeqno{immersed_face_area} circumvents numerical pitfalls: %
While computing $\vx_{0,k}\not\in\vfface_k$, e.g., like in \reffig{face_brackets} (left), poses no challenge in theory, a naive numerical evaluation suffers from severe loss of precision, if the intersection of $\plicintersectionpoints$ and $\vfface_k$ lies far away from the cell, which is possible in fully or non-wetted configurations. In this case, the summands in \refeqn{face_area_sequential_computation} admit values which may differ by several orders of magnitude. %
The computation of the area must thus take into account the intersection topology of the respective face $\vfface_k$ to circumvent this obstacle: it is sufficient to know, if the intersection is outside the convex hull of the polyhedron, to simplify the computation.%
\end{remark}%

\refEqn{immersed_face_area} can be interpreted as follows: if $\secondarysigndist\in[\sdlowertriple{k},\sduppertriple{k}]$, i.e., if $\plicintersectionpoints$ intersects the convex hull of $\vfface_k$, the area is computed from the truncated edges $\cutvfedge_{k,m}$ associated to $\vfface_k$; cf.~the right panel in \reffig{face_truncated_edges}. The negative half-spaces of $\primaryplicplane$ and $\secondaryplicplane$ do not overlap in $\vfface_k$ for $\secondarysigndist<\sdlowertriple{k}$ (left panel in \reffig{face_brackets}), implying that the immersed area $\immersedarea{k}$ can be computed from the original (i.e., non-truncated) edges $\vfedge_{k,m}$. %
If the secondary PLIC plane $\secondaryplicplane$ does not intersect $\cutvfface_k$, one trivially obtains %
\begin{align}%
\immersedarea{k}=%
\begin{cases}%
0&\text{if}\quad\secondarysigndist\leq\sdlowerface{k},\\%
\abs{\vfface_k}-\primaryimmersedarea{k}&\text{if}\quad\secondarysigndist\geq\sdupperface{k}.%
\end{cases}%
\label{eqn:immersed_area_cut}%
\end{align}
\refEqn{immersed_area_cut} induces a status of the faces, which allows to hierarchically infer the topological classification of the polyhedron described in \reffig{degenerate_configurations}. For a given signed distance $t$, the status is %
\textbf{non-wetted} if $t\leq\min_k\sdlowerface{k}$, %
\textbf{fully wetted} if $t\geq\max_k\sdupperface{k}$ and %
\textbf{triple} otherwise.%

The immersed edge length $\ell$ of a non-degenerately intersected edge $\vfedge$ ($\vfstatus[1]{\vfedge}=0$, cf.~\reftab{intersection_status_edge}) with vertices $\vfvert_1\defeq\vx^\vfface_{k,m}$ and $\vfvert_2\defeq\vx^\vfface_{k,m+1}$, renamed for ease of notation, is computed as%
\begin{align}
\ell\fof{\vfedge;\secondarysigndist}=\norm{\vfvert_1-\vfvert_2}%
\begin{cases}
\phantom{1\,-\,}\lengthcoeffconst{}{}+\secondarysigndist\lengthcoefflin{}{}&\text{if}\quad\vfstatus[2]{\vfvert_1}=-1,\\%
{1\,-\,}\lengthcoeffconst{}{}-\secondarysigndist\lengthcoefflin{}{}&\text{if}\quad\vfstatus[2]{\vfvert_1}=1,%
\end{cases}\label{eqn:edge_length_computation}%
\intertext{with the generic coefficients}%
\lengthcoeffconst{}{}=\frac{\iprod{\xbase-\vfvert_1}{\secondaryplicnormal}}{\iprod{\vfvert_2-\vfvert_1}{\secondaryplicnormal}}%
\quad\text{and}\quad%
\lengthcoefflin{}{}=\frac{1}{\iprod{\vfvert_2-\vfvert_1}{\secondaryplicnormal}}.\label{eqn:edge_length_coefficients_generic}%
\end{align}
For non-intersected edges ($\vfstatus[1]{}\neq0$), cf.~\reftab{intersection_status_edge}), one obtains %
\begin{align}
\ell\fof{\vfedge}=%
\begin{cases}
0&\text{if}\quad\vfstatus[1]{\vfedge}\in\set{1,2},\\%
\norm{\vfvert_1-\vfvert_2}&\text{if}\quad\vfstatus[1]{\vfedge}\in\set{-1,-2,3}%
\end{cases}.%
\label{eqn:edge_length_status}%
\end{align}%
\begin{remark}[Regularity]%
The status-based assignment in \refeqn{edge_length_status} corresponds to assigning the left-sided limit to the derivatives of the immersed area $\immersedarea{k}$ with respect to $\secondarysigndist$ for $\secondarysigndist\in\set{\iprod{\vx^\vfface_{k,m}-\xbase}{\secondaryplicnormal}}_{m}$, i.e., when the plane $\secondaryplicplane$ passes through one of the faces' vertices; cf.~subsection~2.2 (\textit{Topological properties of geometric entities}) in \citet{JCP_2021_fbip}. %
\end{remark}
For an edge $\vfedge_{k,m}$ of the original polyhedron, \refeqn{edge_length_coefficients_generic} becomes %
\begin{align}%
\lengthcoeffconst{k}{m}=\frac{\iprod{\xbase-\vx^\vfface_{k,m}}{\secondaryplicnormal}}{\iprod{\vx^\vfface_{k,m+1}-\vx^\vfface_{k,m}}{\secondaryplicnormal}}%
\quad\text{and}\quad%
\nonumber
\\ 
\lengthcoefflin{k}{m}=\frac{1}{\iprod{\vx^\vfface_{k,m+1}-\vx^\vfface_{k,m}}{\secondaryplicnormal}}.\label{eqn:edge_length_coefficients}%
\end{align}
The specific form of \refeqn{edge_length_coefficients_generic} for a truncated edge $\cutvfedge_{k,m}$ (cf.~\refeqn{truncated_edge_assignment} and \reftab{intersection_status_edge} for the respective assignment of $\vfvert_1$ and $\vfvert_2$, depends on the logical status $\vfstatus[1]{}$ (see \reftab{intersection_status_edge}) of its vertices with respect to the \textit{primary} plane $\primaryplicplane$, i.e., %
\begin{equation}
\begin{aligned}
\brackets*{\lengthcoeffconst{}{},\lengthcoefflin{}{}}_{k,m}=%
\begin{cases}
\frac{\brackets*{\iprod{\xbase-\vx^\vfface_{k,m+1}}{\secondaryplicnormal},1}}{\brackets{\beta-1}\iprod{\vx^\vfface_{k,m+1}-\vx^\vfface_{k,m}}{\secondaryplicnormal}}&\text{if}\quad\vfstatus[1]{\vx^{\vfface}_{k,m}}=-1,\\[10pt]%
\frac{\brackets*{\iprod{\xbase-\vx^\vfface_{k,m}}{\secondaryplicnormal},1}}{\beta\iprod{\vx^\vfface_{k,m+1}-\vx^\vfface_{k,m}}{\secondaryplicnormal}}&\text{if}\quad\vfstatus[1]{\vx^{\vfface}_{k,m}}=1,%
\end{cases}%
\end{aligned}
\label{eqn:cutedge_length_coefficients}
\end{equation}
with
\begin{equation*}
\beta=\frac{\primarysigndistref+\iprod{\xbase-\vx^\vfface_{k,m}}{\primaryplicnormal}}{\iprod{\vx^\vfface_{k,m+1}-\vx^\vfface_{k,m}}{\primaryplicnormal}}. %
\end{equation*}
%
%
\begin{remark}[Computational efficiency]%
Note that the coefficients in \refeqn{edge_length_coefficients} and \refeqn{cutedge_length_coefficients} are static in the sense that they are computed only once (after the primary plane $\primaryplicplane$ is positioned), which is advantageous since they are required for all iterations of the root-finding scheme for matching the secondary volume. %
\end{remark}%
So far, we have successfully extended the highly efficient divergence-based volume computation to twice truncated polyhedra, where a given primary and secondary volume is truncated from the polyhedral cell by the sought PLIC planes. The advantage over a simple repetition of the \glqq\textit{single}\grqq{} positioning scheme results from the fact that there is no need to establish the connectivity of the cut polyhedron after the first truncation: the choice of a reference point $\plicbase$ on the contact line leads to a zero contribution of any polygon coplanar to the PLIC planes, hence only polygons on the original polyhedron's faces have to be considered. %
Below, a significant performance gain is shown already for \textsc{Cartesian} meshes in our reference solver. %
As the algorithm is not restricted to cuboidal cells, the gain is expected to be even larger on unstructured meshes involving arbitrary polyhedrons.%
%

%
%
%
%
%
\section{Numerical experiments}\label{sec:numerical_experiments}%
\begin{figure*}[!hbt]%
\null\hfill%
\includegraphics{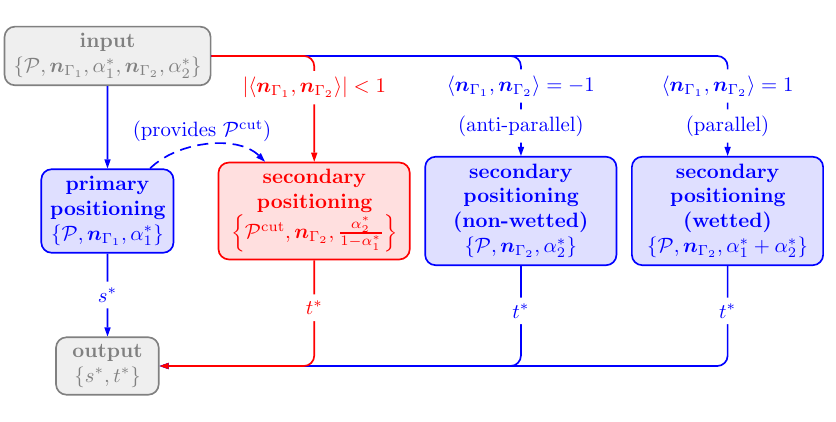}%
\hfill\null%
\caption{Flowchart of the sequential interface positioning algorithm for all three-phase configurations: %
Each of the blue boxes corresponds to performing the positioning of the first truncation and the degenerate configurations by employing the method introduced in \citet{JCP_2021_fbip}, which is applicable to any independent positioning of a single PLIC plane in an arbitrary polyhedron $\protect\polyhedron*$. %
The red box highlights the non-degenerate configuration 
(see~\protect\reffig{degenerate_configurations}), %
which resorts to the adapted volume computation in the cut polyhedron $\cutpolyhedron$ (see~\protect\reffig{sequential_truncation}) for the secondary PLIC positioning developed in the present work.
}%
\label{fig:flowchart}%
\end{figure*}

With the volume computation of a sequentially truncated polyhedron at hand, cf.~\refeqn{volume_fraction_secondary}, \reffig{flowchart} provides a flowchart of the adapted positioning algorithm. %

As stated above, the present paper introduces an extension of the work of \citet{JCP_2021_fbip} in terms of the volume computation for a twice truncated polyhedron with the advantage that no connectivity has to be re-established before the second truncation. A flowchart of the procedure including the non-degenerate as well as degenerate (anti-)parallel cases is shown in \reffig{flowchart}. However, the root-finding part, which was probed with an extensive set of normals and volume fractions for different convex and non-convex polyhedra, remains unchanged in comparison to \cite{JCP_2021_fbip}. In terms of the number of polyhedron truncations, constituting the relevant measure for computational effort, the root-finding algorithm of \citet{JCP_2021_fbip} exhibits high efficiency: on average, the positioning requires \num{1}--\num{2} polyhedron truncations, thus outperforming existing methods; cf.~\refapp{implicit_bracketing}. Hence, we consider it sufficient to conduct only a small set of numerical experiments. %
As pointed out by \citet{XXX_2020_ivof}, a meaningful assessment requires to consider the limiting cases of almost empty ($\refvof\approx0)$ or full ($\refvof\approx1$) cells, where, for two planes to be positioned, all combinatorial configurations need to be examined. %
Furthermore, while the volume fraction is invariant under linear transformations, the normal $\plicnormal$ becomes distorted in general. %
In numerical experiments, this effect must be compensated by a sufficiently resolved sample set for $\plicnormal$; cf.~\cite[note~1.2]{JCP_2021_fbip}. %
Subsection~\ref{subsec:experiment_design} covers the details of artificially generated instances. %
In order to assess the performance of our algorithm in comparison to an existing method, we restrict ourselves to cubes for the numerical tests. The reason for this is twofold: %
\begin{enumerate}
\item The root-finding method employed was shown to be highly efficient for several classes of polyhedra \cite{JCP_2021_fbip}. Thus, only the augmented volume computation, which allows two subsequent positioning tasks without re-establishing the connectivity, remains to be assessed. Due to the face-based formulation, the computational demand of the volume computation increases linearly with the number of faces, but is otherwise independent of the polyhedron topology. Hence, it is sufficient to consider a fairly simple polyhedron. %
\item The flow solver, within which the reference algorithm was originally implemented, resorts to a \textsc{Cartesian} mesh of cuboidal cells. The existing scheme, which will serve as a performance reference, combines a polyhedron decomposition for the volume computation with an accelerated bisection for the root-finding, a method only applicable to cuboids; cf.~section~\ref{sec:decomposition_approach} and \refapp{decomposition_approach} for a detailed description. %
\end{enumerate}
\subsection{Design of the numerical experiments}\label{subsec:experiment_design}%
\begin{figure}[!tb]%
\null\hfill%
\includegraphics{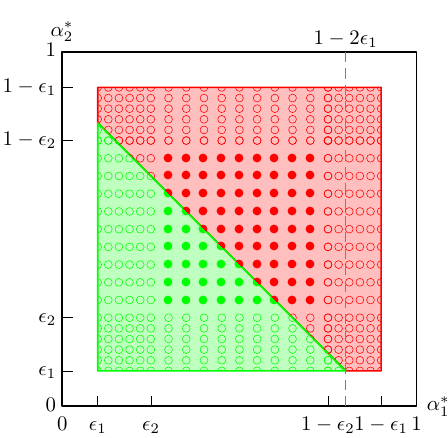}%
\hfill\null%
\caption{Illustration of sampling strategy for volume fractions $(\primaryrefvof,\secondaryrefvof)\in{[0,1]^2}$ for tolerances $\epsilon_1,\epsilon_2$ with admissible (log: \textcolor{green}{$\circ$}, lin: \textcolor{green}{$\bullet$}) and non-admissible (log: \textcolor{red}{$\circ$}, lin: \textcolor{red}{$\bullet$}) configurations. %
Note that both volume fractions refer to the volume of the original polyhedron $\protect\polyhedron*$. %
}%
\label{fig:experiment_design}%
\end{figure}

In the three-material case, the design of a meaningful numerical experiment, i.e., the sample sets for the normals and volume fractions, requires some preliminary considerations. %
While one can choose the normals $\primaryplicnormal$ and $\secondaryplicnormal$ independently, the volume fractions must be chosen carefully, based on two user-defined tolerances $0<\epsilon_1,\epsilon_2\ll1$ with $\epsilon_1>\epsilon_2$. %
With $\epsilon_1$ being the most restrictive tolerance, \refnote{floating_point_operations} induces the constraints %
\begin{equation}
\begin{aligned}
\epsilon_1 & \leq\primaryrefvof\leq1-\epsilon_1,\quad \\%
\epsilon_1 & \leq\secondaryrefvof\leq1-\epsilon_1 %
\quad\text{and}\quad \\%
\epsilon_1 & \leq\primaryrefvof+\secondaryrefvof\leq1-\epsilon_1%
.%
\end{aligned}
\label{eqn:vof_numerical_conditions}
\end{equation}
All volume fraction combinations of the two volume fractions $\primaryrefvof$ and $\secondaryrefvof$ are indicated by $\circ$ and $\bullet$ in \reffig{experiment_design}, but some physical constraints, tolerances and considerations are still to be added for the numerical assessment of three phase interface reconstruction:
In the assessment of the two-phase case, corresponding to a single plane to be positioned, \citet{JCP_2021_fbip} investigated volume fractions $\refvof$ from a sample set $\vofset{}\subset[\epsilon_1,1-\epsilon_1]$, which admits logarithmic spacing indicated by $\circ$ for volume fractions in the vicinity of zero and one, i.e., for $\refvof\in{[\epsilon_1,\epsilon_2]\cap[1-\epsilon_2,1-\epsilon_1]}$, and linear spacing indicated by $\bullet$ otherwise, i.e., in ${[\epsilon_2,1-\epsilon_2]}$. Those limits providing a threshold for the detection of a phase are employed for both volume fractions in this work for classifying full, empty and three-phase cells. %
The tolerance $\epsilon_1$ numerically replaces zero and $\epsilon_2$ indicates the transition from logarithmic to linear spacing. %
As in \cite{JCP_2021_fbip}, for the remainder of this manuscript, we employ $\epsilon_1=\num{e-9}$ and $\epsilon_2=\num{e-5}$. %
A naive extension for $(\primaryrefvof,\secondaryrefvof)$, however, comprises configurations violating the condition $\primaryrefvof+\secondaryrefvof\leq1-\epsilon_1$, i.e., overfull cells or cells, which are numerically detected as two phase cells are not allowed. %
Besides the admissibility, \reffig{experiment_design} illustrates another distinct feature of the sequential positioning: due to the numerical zero-detection discussed in \refnote{floating_point_operations}, the accuracy of the reconstruction looses its symmetry. %
This numerical artefact is best illustrated by considering the case $\primaryrefvof=1-\epsilon_1$. Compressing the second and third expression in \refeqn{vof_numerical_conditions} produces the conflicting condition $\epsilon_1\leq\secondaryrefvof\leq0$, implying that, in a numerical sense, the secondary phase will be neglected because the cell is fully occupied by the primary phase. In order to enforce a tolerance of $\epsilon_1$ for the secondary phase, the upper bound for the primary volume fraction \textit{effectively} reduces to $1-2\epsilon_1$, where an analogous statement holds vice versa. Thus, \refeqn{vof_numerical_conditions} becomes %
\begin{equation}
\begin{aligned}
\epsilon_1 & \leq\primaryrefvof\leq1-2\epsilon_1,\quad \\%
\epsilon_1 & \leq\secondaryrefvof\leq1-2\epsilon_1%
\quad\text{and}\quad \\%
\epsilon_1 & \leq\primaryrefvof+\secondaryrefvof\leq1-\epsilon_1%
, 
\end{aligned}
\label{eqn:vof_numerical_conditions_effective}%
\end{equation}
corresponding to the shaded green region in \reffig{experiment_design}. %
With the sample sets for the normals %
\begin{align*}
\normalset{& M_{\plicnormal}}\defeq\biggl\{\brackets*[s]{\cos\varphi\sin\theta,\sin\varphi\sin\theta,\cos\theta}\transpose: \\
 & \brackets{\varphi,\theta}\in\frac{\pi}{2M_{\plicnormal}}\brackets[s]{1,2,\dots,2M_{\plicnormal}}\times\frac{\pi}{M_{\plicnormal}}\brackets[s]{0,1,\dots,M_{\plicnormal}}\biggr\},%
\end{align*}
which is parametrized by the spherical angles $\theta$ and $\varphi$, and volume fractions 
\begin{equation*}
\begin{aligned}
&\vofset{M_{\refvof}}   \defeq \\
& \set*{10^{-4}+\frac{m-1}{M_{\refvof}-1}\brackets*{1-2\cdot10^{-4}}:1\leq m\leq M_{\refvof}}\cup\\%
& \set*{10^{-9},  10^{-8},10^{-7},10^{-6},10^{-5}}\cup \\%
&  \set*{1-10^{-5},  1-10^{-6},1-10^{-7},1-10^{-8}, 1-2\cdot10^{-9}}\text{,}%
 \end{aligned}
\end{equation*}
one obtains the following sets for the tupels of normals and volume fractions: %
\begin{equation}
\begin{aligned}
\mathcal{S}^2_{\refvof} & \defeq\set{(x,y)\in\vofset{}\times\vofset{}:x+y\leq1-\epsilon_1}%
\quad\text{and}\quad \\%
\mathcal{S}^2_{\plicnormal} & \defeq\normalset{}\times\normalset{}.%
\end{aligned}
\end{equation}
Note that the size\footnote{At the poles, i.e., for $\theta\in\set{0,\pi}$, the azimuthal angle $\varphi$ carries no information. Hence, the size of $\normalset{M_{\plicnormal}}$ effectively reduces from $2M_{\plicnormal}\brackets{M_{\plicnormal}+1}$ to $2+2M_{\plicnormal}\brackets{M_{\plicnormal}-1}$.} of their dyadic product, corresponding to the total number of instances, grows very quickly: %
choosing the moderate resolutions $M_{\plicnormal}=10$ and $M_{\refvof}=20$, one already obtains $\abs{\mathcal{S}^2_{\refvof}\times\mathcal{S}^2_{\plicnormal}}\approx\num{1.45e7}$ samples. %
\subsection{Baseline for the performance comparison: A decomposition-based accelerated bisection approach}\label{sec:decomposition_approach}%
The present method is compared to an existing method in the multi-phase flow solver FS3D based on polyhedron decomposition for the volume computation and an accelerated bisection for the root-finding; cf.~\refapp{decomposition_approach} and \citet{Potyka_2018} (in German) for details. In a nutshell, the approach can be summarized as follows:%
\begin{enumerate}
\item Rotate the coordinate system such that $\iprod{\plicnormal}{\ve_1}\geq\iprod{\plicnormal}{\ve_2}\geq\iprod{\plicnormal}{\ve_3}$. %
\item Exploit congruency, i.e., $\set{\refvof,\plicnormal}\equiv\set{1-\refvof,-\plicnormal}$ for $\refvof<\frac{1}{2}$
\item The volume function can be explicitly inverted using the formulae from \refapp{cube_explicit_volume_inverse} for a cuboid cell. %
This implies that the tentative positioning of the primary plane in all cases and of the secondary plane in the (anti-)parallel degenerate cases (blue boxes in \reffig{flowchart_decomposition}) is carried out directly with the explicit inversion. The primary position $\primarysigndistref$ is found with this explicit method described in detail in \refapp{cube_explicit_volume_inverse}. %
\item The secondary position $\secondarysigndistref$ is found iteratively inside the residual polyhedron with an accelerated bisection which makes use of the two-phase positioning of an extended volume fraction inside the original polyhedron $\polyhedron*$ and an error calculation of the position $\secondarysigndist$ inside the remaining polyhedron $\cutpolyhedron$, cf.~\refapp{decomposition_approach}. The volume error constitutes a lower limit for the change required for the iteration limits in the following iteration. Thus, utilizing this lower bound to update the iteration boundaries accelerates the bisection root-finding.
\item The computation of the truncated volume $\polyvof_2\fof{\secondarysigndist;\cutpolyhedron}$ relies on tetrahedral decomposition of convex polyhedra.%
\end{enumerate}
Both the volume computation and the intersection of a polyhedron with a plane require establishing the full connectivity (i.e., vertex to edge, edge to face, face to cut polyhedron), corresponding to multiple computationally expensive comparing and sorting operations required in the accelerated bisection approach, but not in our new algorithm. While the algorithm, we compare the new method to, may be feasible for the cuboid meshes the algorithm was designed for, the feasibility quickly degrades for non-convex meshes. %

\subsection{Numerical results}\label{subsec:numerical_results}%
\begin{figure*}[!tb]%
\centering%
\null\hfill%
\includegraphics[page=1]{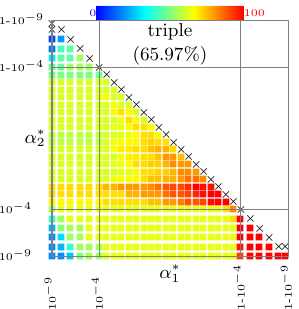}%
\hfill%
\includegraphics[page=2]{\thisfigurepath coverage}%
\hfill%
\includegraphics[page=3]{\thisfigurepath coverage}%
\hfill\null%
\caption{Relative incidence of intersection topology (in \%; cf.~\protect\reffig{degenerate_configurations}) as function of primary and secondary volume fraction. E.g., most of the non-wetted configurations (rightmost) are obtained if both the primary and secondary volume fraction are very small, i.e., $\num{e-9}\leq\primaryrefvof,\secondaryrefvof\leq\num{e-4}$.}%
\label{fig:experiment_topo_coverage}%
\end{figure*}

\begin{figure*}[!tb]%
\centering
{%
\null\hfill%
\includegraphics[page=1]{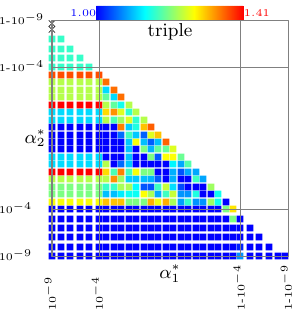}%
\hfill%
\includegraphics[page=2]{\thisfigurepath iterations}%
\hfill%
\includegraphics[page=3]{\thisfigurepath iterations}%
\hfill\null%
}%
\\%
{%
\null\hfill%
\includegraphics[page=4]{\thisfigurepath iterations}%
\hfill%
\includegraphics[page=5]{\thisfigurepath iterations}%
\hfill%
\includegraphics[page=6]{\thisfigurepath iterations}%
\hfill\null%
}%
\caption{Number of polyhedron truncations $N_{\mathrm{av}}$ (see~\refeqnsafe{result_average_definition}) required for positioning the secondary plane $\protect\secondaryplicplane$ as function of primary and secondary volume fraction for present divergence-based (top row) and decomposition-based positioning (bottom row), disaggregated by intersection topology; cf.~\reffig{degenerate_configurations}. Note the non-coinciding ranges, respectively indicated to the left and right of the color bars. Note the significant improvement (on average, 1 to 1.41 iterations are required, as compared 4.59 to 31.97) achieved by the proposed algorithm for the triple-configuration, which is the most common one in the stand-alone assessment of the old and new algorithms.}%
\label{fig:experiment_result}%
\end{figure*}

In what follows, we conduct numerical experiments based on the data sampling strategy outlined in subsection~\ref{subsec:experiment_design}. %
First, recall from \reffig{flowchart} that the sequential character of the algorithm implies that the primary positioning is carried out independently. Hence, the number of secondary polyhedron truncations, i.e., the number of truncations applied to $\cutpolyhedron$, constitutes the relevant performance measure. %
%
We aggregate the number of secondary truncations $N$ for each pair of volume fractions $\brackets{\primaryrefvof,\secondaryrefvof}$ by defining the average over the normals, i.e., %
\begin{align}
N_{\mathrm{av}}\fof{\primaryrefvof,\secondaryrefvof}\defeq%
\frac{1}{\abs{\mathcal{S}^2_{\plicnormal}}}\sum\limits_{i=1}^{M_{\plicnormal}}\sum\limits_{j=1}^{M_{\plicnormal}}
{N\fof{\primaryrefvof,\secondaryrefvof,%
\vn_{\plicplane,i},\vn_{\plicplane,j}%
}}.%
\label{eqn:result_average_definition}%
\end{align}
For the interpretation of the results, it is instructive to disaggregate the averaged number of truncations based on the topological classification given in \reffig{degenerate_configurations}, which can be one of the following: \textit{triple}, \textit{fully wetted} or \textit{non-wetted}. %
Recall from \reffig{flowchart} that the (anti-)parallel degenerated configurations are treated independently from the non-degenerate cases. %
\begin{note}[Performance impact of input data in practical applications]\label{note:application_parallel}%
It is worth noting that the relative incidence of the intersection topologies strongly depends on the input data, i.e., on the volume fractions and normals; cf.~\reffig{experiment_topo_coverage}. %
The physical properties of the involved fluid phases induce individual and highly dynamic spectra of configurations: e.g., consider a droplet (secondary phase) impinging on a structured super-hydrophilic wall (primary phase) in ambient gas: right before the impingement, there are only non-wetted (cyan) configurations. Once the contact line starts to form, the majority of the three-phase cells contains a triple line (green). With the wetted area and the number of fully wetted cells (blue) being of codimension one and the advancing contact line being of codimension two, the ratio quickly shifts towards fully wetted cells. 
\\
\begin{center}%
\includegraphics[page=1, scale=1.1]{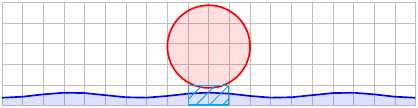}%
\\
\includegraphics[page=2, scale=1.1]{\thisfigurepath application_illustration}%
\end{center}%
In view of an application within a parallelized solver, where the phase boundaries and triple line may comprise several processes, load balancing becomes relevant. Thus, it is of paramount importance to ensure a robust performance with respect to the intersection topology. %
\end{note}
\refFig{experiment_topo_coverage} depicts the relative incidence of topological classifications encountered, as a function of the volume fractions $(\primaryrefvof,\secondaryrefvof)$, where "$\times$" denotes zero incidence, i.e.\ for each combination of $(\primaryrefvof,\secondaryrefvof)$, the values in the panels add up to 100 (\%). Also, recall from \refeqn{vof_numerical_conditions_effective} that the choice of volume fractions is restricted, reproducing the lower triangluar structure of \reffig{experiment_design}. %
One obtains the accumulated ratio of incidences (\% in parentheses) by further averaging over all combinations of volume fractions. The accumulated ratio of incidences is the fraction of the respective classification with respect to the total number of instances $\abs{\mathcal{S}^2_{\refvof}\times\mathcal{S}^2_{\plicnormal}}$. As visible in the leftmost panel, roughly two thirds (65.97\%) of all instances under consideration feature a triple line. %
No \textit{fully wetted} configurations were encountered (center panel in \reffig{experiment_topo_coverage}) for instances with $\secondaryrefvof\lessapprox\num{e-4}$. %
Analogously, \textit{non-wetted} configurations do not occur for $\secondaryrefvof\gtrapprox1-\num{e-4}$ (right panel in \reffig{experiment_topo_coverage}). %
It is crucial to note that the depicted relative incidence strongly depends on the underlying sample set of the normals $\mathcal{S}^2_{\plicnormal}$, given in subsection~\ref{subsec:experiment_design}. %
Also, \reffig{experiment_topo_coverage} exemplifies the performance considerations of \refnote{application_parallel}: while triple configurations can be found for nearly all combinatiosn of $\primaryrefvof$ and $\secondaryrefvof$, a wetted configuration requires the secondary volume fraction to exceed a certain threshold. %
\refFig{experiment_result} compares the average number of secondary truncations $N_{\mathrm{av}}$ (see~\refeqn{result_average_definition}) of the proposed algorithm to the existing decomposition-based accelerated bisection scheme of \citet{Potyka_2018}, cf.~\refapp{decomposition_approach}. %

The main implications can be cast as follows: %
\begin{enumerate}
\item The proposed method maintains the high efficiency of the single PLIC plane positioning of \cite{JCP_2021_fbip}. 
This becomes evident by comparing the maximum number of truncations, indicated to the left and right of the respective color bar. This maintained high efficiency of the root-finding is to be expected, since the divergence-based root-finding approach is the same while the volume computation was enhanced for sequential truncation. Our approach requires, on average, between \num{1} and \num{2} truncations for the positioning of the secondary plane for all possible intersection classifications depicted in \reffig{degenerate_configurations}. Note that, with \num{1} to \num{1.41} truncations on average (depending on the volume fraction combination), our algorithm performs best in the presence of a triple-line, accounting for roughly two thirds of the instances. %
\item The average number of truncations slightly increases for wetted intersection configurations as the secondary volume fraction $\secondaryrefvof$ approaches the value below which the intersection topology changes ($\times$ in \reffig{experiment_result}). An analogous phenomenon occurs for non-wetted intersections when the sum of volume fraction approaches $1-\epsilon_1$, i.e., for almost full cells, as well as for two distinct combinations ($\set{1-\num{e-4},\num{e-9}}$ and $\set{1-\num{e-4},\num{e-8}}$). %
Nevertheless, the differences are significantly smaller than for the algorithm we compare to. It highlights the robustness of the proposed efficient sequential PLIC algorithm: besides the aforementioned the performance admits no evident pattern for variation of (combinations of) $\primaryrefvof$ and $\secondaryrefvof$. %
In three-phase cells, the execution time admits an about even distribution, which is beneficial in terms of load balance within a parallelized application. %
\item The number of truncations required by the accelerated decomposition-based bisection approach\footnote{See~\refapp{decomposition_approach} for a full description.} for configurations of \textit{triple}-type shown in \reffig{experiment_result} (bottom row, left) strongly depends on the combination of volume fractions. The explanations for this strong dependency are the following: %
\begin{enumerate}
\item If the volume fraction $\polyvof_2\fof{\secondarysigndist;\polyhedron*}$ and the calculated error $\Delta\polyvof_2\fof{\secondarysigndist;\cutpolyhedron}$ differ by more than an order of magnitude, the effect of $\Delta\polyvof_2$ on the update of $\polyvof_2\fof{\secondarysigndist;\polyhedron*}$ diminishes. In such cases, the algorithm degenerates to pure bisection for $\polyvof_2<\polyvof_1$, as the update with the error has virtually no effect on the change of the iterated volume fraction inside the whole cell.
\item If $\primaryrefvof$ is large compared to $\secondaryrefvof$, the initial bracketing interval for the iteration is quite large, while the interval is small at small $\primaryrefvof$, cf.~\reffig{face_bracket_decomposition}.
\item The convergence is slow when the sought value is relatively close to an initial bracketing value, thus the opposite bound requires multiple updates in the beginning.%
\item A superposition of aforementioned negative effects degrades the convergence for some cases. %
\end{enumerate}
\end{enumerate}

\begin{note}[Topological classification]%
\label{note:TopologicalClassification}
The figure* below shows the number of intersection topology mismatches between the present and decomposition-based approach (in \%, relative to the number of possible normal combinations $\abs{\mathcal{S}^2_{\plicnormal}}$, where $\times$ indicates zero mismatches; cf.~\refeqn{result_average_definition}) as function of primary and secondary volume fraction $(\primaryrefvof,\secondaryrefvof)$ (overall quota: 0.47\%). Note that the absolute accuracy of the position is not degraded by a mismatch.\\%
\null\hfill
\includegraphics[scale=1.2]{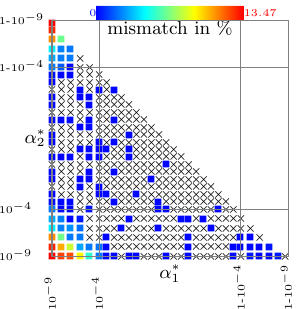}
\hfill\null%
\\
\null\hfill
\includegraphics{\thisfigurepath mismatch_illustration}
\hfill\null%
\\%
While the topological classification of the present algorithm is based on the tolerance for the signed distance $\secondarysigndist$ ($\zerotol=\num{e-14}$), cf.~\refnote{volumepositiontol}, the decomposition-based approach employs the tolerance used for the volume fraction ($\epsilon=\num{e-12}$): if the error $\Delta\polyvof_2$ in the volume fraction for an assumed \textit{non-wetted} ($\secondarysigndist^{\prime}$) or \textit{wetted} ($\secondarysigndist^{\prime\prime}$) configuration falls below $\epsilon=\num{e-12}$, the respective classicfication is assigned. %
In all other cases, i.e., if both $\abs{\Delta\polyvof_2\fof{\secondarysigndist^{\prime}}}$ and $\abs{\Delta\polyvof_2\fof{\secondarysigndist^{\prime\prime}}}$ exceed $\epsilon$, the \textit{triple} classification is at hand; cf.~\reffig{flowchart_decomposition} and \refapp{decomposition_approach}. %
This conceptual difference may produce dissenting classifications in the vicinity of a transition between two classifications. Note that, in general, no fixed choice of length and volume tolerances capable of resolving this artefact exists. %
From the right panel in the above figure*, the rationale emerges from a simple geometric consideration: for some given signed distance, say $\bar{\secondarysigndist}$, a small variation $\delta\secondarysigndist$ translates to a change of the volume fraction $\delta\polyvof=\polyvof\fof{\bar{\secondarysigndist}+\delta\secondarysigndist}-\polyvof\fof{\bar{\secondarysigndist}}\approx\delta\secondarysigndist\frac{\abs{\polyhedron*\cap\secondaryplicplane\fof{\bar{\secondarysigndist}}}}{\abs{\polyhedron*}}$. %
In other words, for a consistent choice, the length and volume tolerances would have to be adjusted to the dynamic amplification $\nicefrac{\delta\polyvof}{\delta\secondarysigndist}$.
\end{note}

\subsection{Embedded performance assessment}\label{Subsec:PerfInFS3D}%
\newcommand{\perfplotwidth}{6cm}%
The software FS3D was originally developed by \citet{PhD_2004_Rieber} for the Direct Numerical Simulation (DNS) of incompressible two-phase flows with sharp interfaces. %
Among the applications of FS3D are, e.g., freely rising droplets \cite{JFM_2015_dbob}, binary droplet collision~\cite{PF_2012_dnso,CIT_2012_dnsb,JFM_2016_nsoh,AM_2019_ttps}, falling films~\cite{IJMF_2012_iost} and thermocapillary flows~\cite{IJMF_2015_dnso}; the interested reader may also consider the comprehensive overview of \citet{EISENSCHMIDT2016}. %
FS3D employs a time-explicit finite volume approach on a \textsc{Cartesian} mesh with a staggered arrangement of variables \cite{PF_1965_ncot} to numerically solve a one-field formulation of the \textsc{Navier-Stokes} equations. %
The interface is captured by the VoF approach of \citet{HIRT1981}, where the volume fraction $\alpha_i(t,\vx)$ encodes the belonging of a spatial point $\vx$ to the phases $\Omega_i^\pm(t)$. The numerical treatment of rigid bodies was incorporated by \citet{RAUSCHENBERGER2015} in terms of an \textsc{Eulerian} framework. %
Its dimensionally split transport \cite{STRANG1968} resorts to the PLIC method of \citet{Youngs1982}. An extension of the PLIC normal orientation algorithm to three phases similar to the well-known LVIRA method with additionally accounting for some peculiarities close to the contact line was introduced into FS3D based on the method of \citet{JCP_2016_atdv}. This three-phase normal orientation strategy was combined with the here presented PLIC positioning method in the multi-phase simulation code FS3D. %
The cuboidal domain $\Omega=\Omega^+\fof{t}\cup\Omega^-\fof{t}\cup\iface\fof{t}$ is decomposed into $N_{i}\times N_{j}\times N_{k}$ fixed control cubic volumes $\Omega_{ijk}$ of size $\Delta x^3$, which are equally distributed among $P_{i}\times P_{j}\times P_{k}$ processors. Additional ghost control volumes at the interior (between processors) and exterior domain boundaries allow for a convenient integration of the boundary conditions into the discretization schemes. \textbf{M}essage \textbf{P}assing \textbf{I}nterface (MPI) is used for the exchange of information between processors, each associated with an MPI process.
While this strategy greatly advantages the iterative multi-grid approach applied to solve the pressure projection, it may result in an uneven distribution of interface cells, which translates to load imbalances among the processes. As the three-phase contact line is of codimension two in a three-dimensional space, the number of processes containing three-phase cells does not scale with increasing even decomposition among more MPI-processes. Although also the three-phase cells are distributed among more processes, then even more processes have to wait for the completion of the three-phase positioning. Our new fast algorithm ensures that this idling time becomes insignificant compared to the overall runtime. Along with the runtime also the power consumption of the processors for a three-phase simulation reduces significantly. %
\begin{figure}[!tb]%
\def\figwidth{7cm}%
\centering
\includegraphics[width=\figwidth]{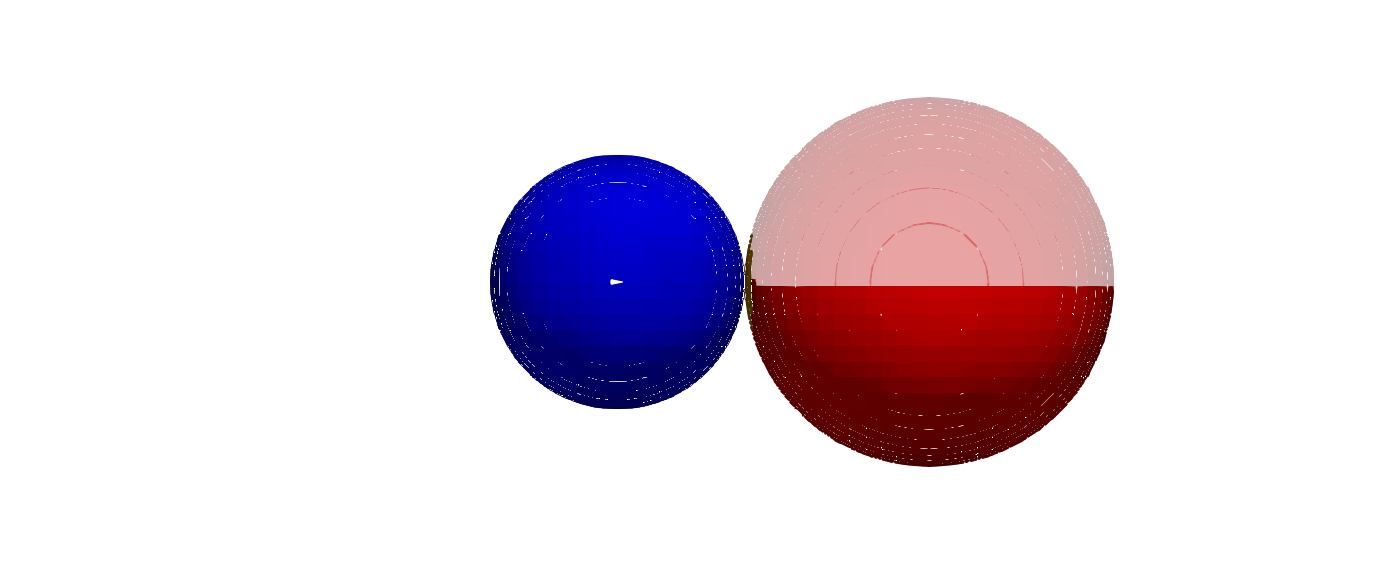}%
\\
\makebox[\figwidth][c]{time step 0}%
\\
\includegraphics[width=\figwidth]{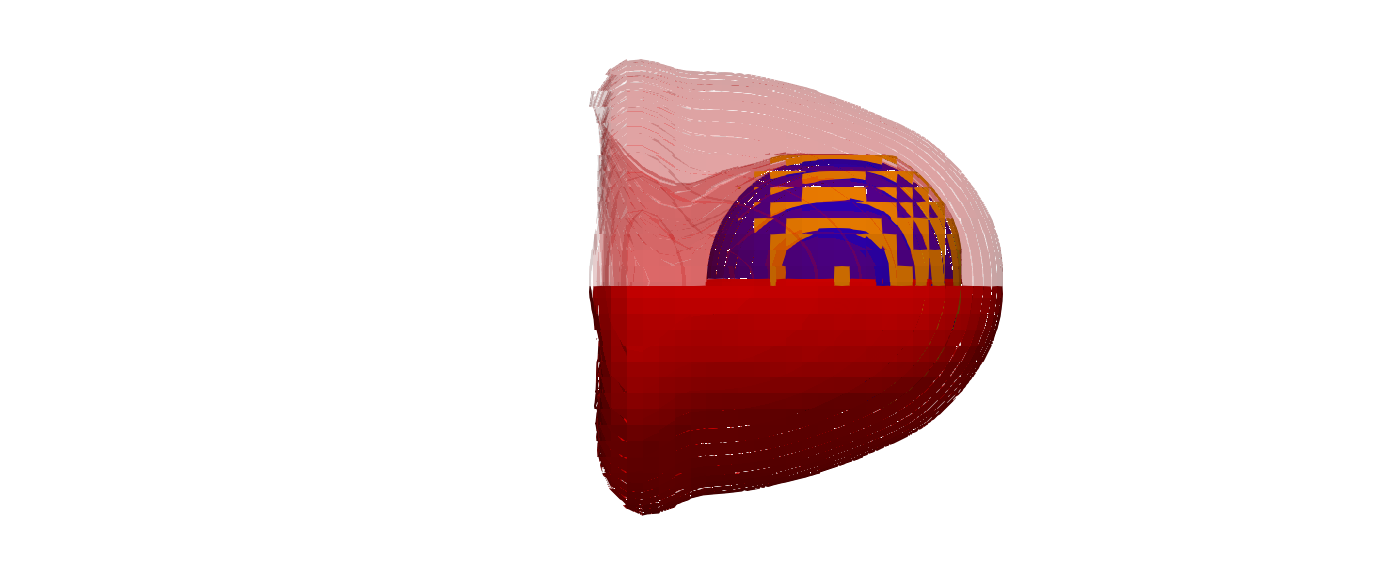}%
\\
\makebox[\figwidth][c]{time step 300}%
\\
\includegraphics[width=\figwidth]{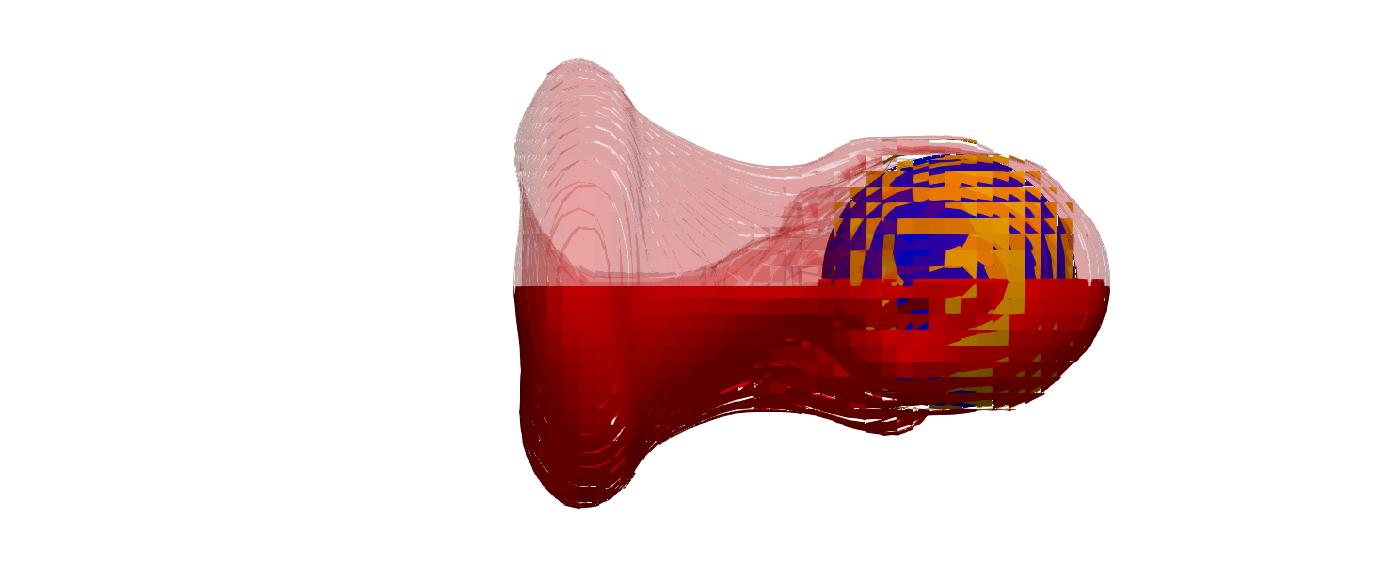}%
\\
\makebox[\figwidth][c]{time step 900}
\\
\includegraphics[width=\figwidth]{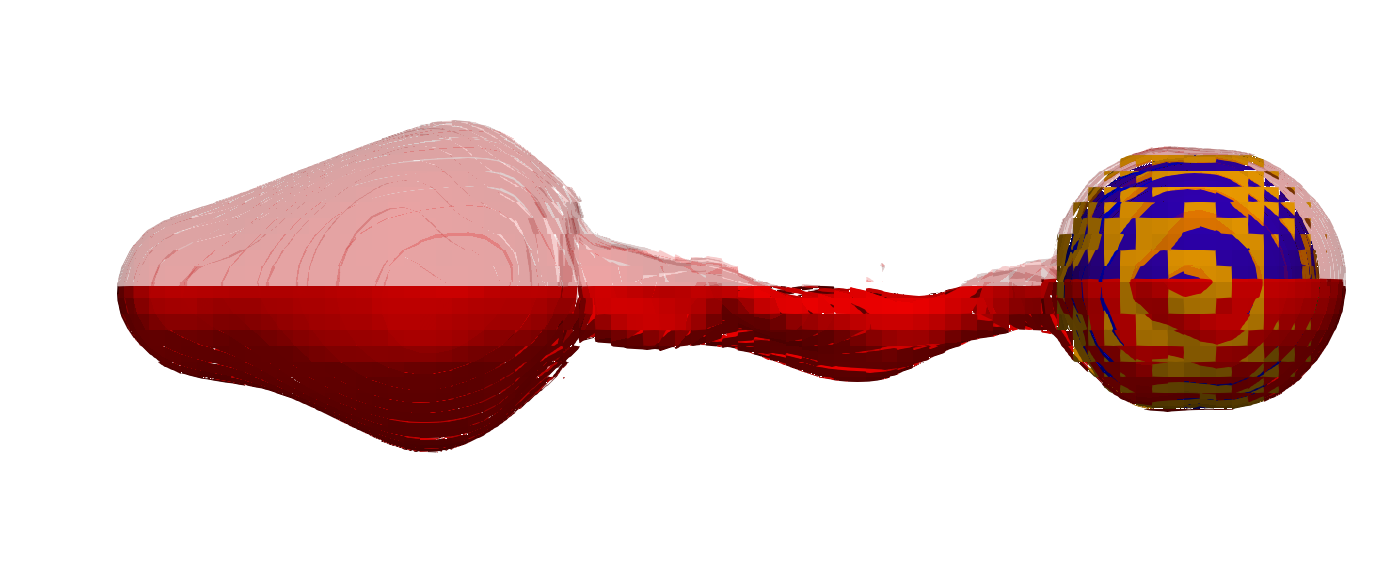}%
\\
\makebox[\figwidth][c]{time step 1800}%
\caption{Visualization of the results for 128 MPI processes: exhibiting air entrapment, large deformations and thin films. Blue PLIC patches represent the solid sphere, whereas the liquid-air interface is colored in red for two-phase cells and orange for three-phase cells. The upper part of the fully three dimensional simulation results is shown transparent.}%
\label{fig:mpiprocs_vis}%
\end{figure}
In order to examine the performance of the proposed positioning algorithm in a prototypical simulation, we consider a collision of a solid sphere and a liquid droplet in ambient gas; see appendix \ref{app:FS3Dcollisionsetup} for further details. 
We vary the equidistant spatial resolution $N$ along with the number of MPI processes $P$ per spatial direction for a fixed setup such that $\nicefrac{N_i}{P_i}=\nicefrac{N_j}{P_j}=\nicefrac{N_k}{P_k}=32$. This ensures that each process contains $32^3$ cells, which reduces the performance impact of the remaining components of the flow solver FS3D. %

\begin{note}[Computational infrastructure]%
\label{note:CompInfrastructure}
The supercomputer Hawk at HLRS in Stuttgart, on which the test cases were run, offers nodes with a maximum of 128 cores. In order order to eliminate the impact of node-to-node communication, the test cases in \reftab{PerformancePLICinFS3D} were configured to fit on a single node. Note that (for the utilized machine) $128$ processes already exhaust the bandwidth of intra-node memory operations. Thus no linear scaling can be expected as the number of processes is increased. %
However, a comparison between the numbers of processors (corresponding each to an MPI process) is irrelevant for the following assessment and has no impact on the results: A comparison of the two methods (in the otherwise identical simulation setup) is the goal of the embedded assessment within an existing flow solver.%
\end{note}%
The strong deformation (cf.~\reffig{mpiprocs_vis}) provides a prototypical spectrum of interface configurations with a representative distribution among the processes, which nicely serves the purpose of demonstrating the performance benefits induced by embedding the proposed algorithm in a parallelized flow solver. %
Note that the spatial resolutions employed in this assessment are not sufficient for full spatial resolution of the smallest process scales. %
Hence, the simulation results do not show a comprehensive picture of the physical collision process, but serve the purpose of demonstrating the reconstruction's performance well. %
\begin{remark}[Performance measurement]%
\label{remark:performanceMeassurment}
In order to limit the influence of the hardware and the operating system, the comparison of the code with the previous and proposed version of the algorithm was performed subsequently on the same node. However, small perturbing influences, e.g., interventions of the operating system or changes in the processor frequency, are possible sources for noise and performance deviations within single time-steps in the measurements. %
\end{remark}%
\begin{figure*}[!tb]%
\setlength{\unitlength}{\linewidth}
\begin{picture}(1,0.85)
\put(-0.01,0.5){\includegraphics[page=1]{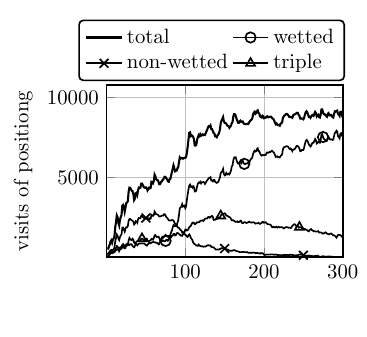}}

\put(0.346,0.5){\includegraphics[page=2]{\thisfigurepath performance_balance_3phcelldistribution.pdf}}

\put(0.664,0.5){\includegraphics[page=3]{\thisfigurepath performance_balance_3phcelldistribution.pdf}}

\put(0.025,0.25){\includegraphics[page=1]{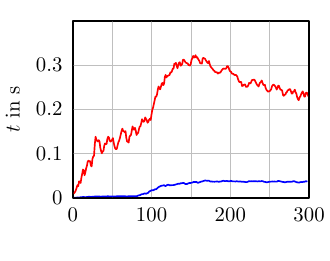}}%

\put(0.346,0.25){\includegraphics[page=2]{\thisfigurepath performance_plic_new}}%

\put(0.675,0.25){\includegraphics[page=3]{\thisfigurepath performance_plic_new}}%

\put(-0.01,0.0){\includegraphics[page=1]{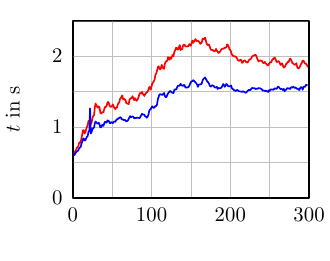}}%
\put(0.346,0.0){\includegraphics[page=2]{\thisfigurepath performance_cycle_new}}%
\put(0.664,0.0){\includegraphics[page=3]{\thisfigurepath performance_cycle_new}}%

\end{picture}
\caption{Execution time (\thisisrevision{top: three-phase cell distribution, middle: PLIC reconstruction, bottom: entire time-step}) summed over all processes (left to right: $2,16$ and $128$ MPI processes) as a function of number of time steps of the flow solver using the previous (red) and proposed (blue) reconstruction algorithm. Note that the CFL condition induces a resolution-dependent number of simulated time steps for ${[0,t_{\mathrm{max}}]}$ with fixed $t_{\mathrm{max}}$. \thisisrevision{The number of calls of the positioning routine for each of the three-phase configurations is shown in the bottom row. The performance measurement of the positioning of the old version mainly follows the number of three-phase cells with the triple configuration. The performance measurement of the new version follows the total number of three-phase cells. This confirms the findings of the stand-alone assessment of the performance in sec.~\ref{subsec:numerical_results}.}}%
\label{fig:runtime}%
\end{figure*}
\begin{table*}[!tb]%
\caption{Comparison of the computational speed-up and average share of reconstruction time for the previous and proposed algorithm embedded in FS3D. The speed-up includes idle times of processes without three-phase cells while the average percentage of reconstruction compared to the total runtime only includes the reconstruction runtime (henceforth called \textit{effective runtime}). While the effective runtime decreases significantly with the resolution, the speed-up stays roughly constant due to the uneven distribution of three-phase cells among the MPI processes. The new algorithm reduces the effective reconstruction runtime by an order of magnitude, thus also the load imbalance becomes insignificant. \thisisrevision{See \refnote{CompInfrastructure} and Remark~\ref{remark:performanceMeassurment} for details on the computation infrastructure used.}}%
\label{tab:PerformancePLICinFS3D}%
\centering%
\renewcommand{\arraystretch}{1.25}%
\begin{tabular}{cc|c c|c c}%
\multicolumn{6}{c}{}\\%
\hline
\textbf{processors}  & \textbf{resolution}  & \multicolumn{2}{c|}{\textbf{speed-up}} & \multicolumn{2}{c}{\textbf{average}$\brackets*{\frac{t_{\mathrm{reconstruct}}}{t_{\mathrm{total}}}}$} \\%
&& reconstruction & total & previous & proposed \\%
\hline
2&  64$\times$32$\times$32   & 89\% & 21\% & 12.6\%  & 1.8 \% \\
16& 128$\times$64$\times$64   & 91\% & 17\% &  6.4\%  & 0.67\% \\
128& 256$\times$128$\times$128 & 95\% & 18\% &  0.90\% & 0.06\% \\
\end{tabular}
\end{table*}
\begin{figure}[!tb]%
	\centering
\includegraphics[page=1]{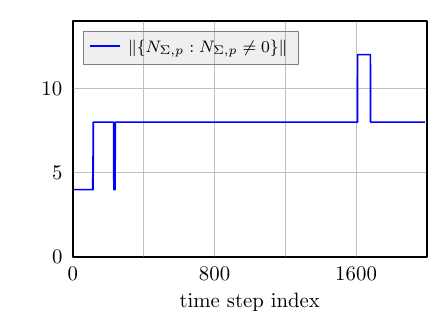}%
\\
\includegraphics[page=2]{\thisfigurepath performance_balance_new}%
\caption{Number of processes containing a non-zero number of three-phase cells $N_{\Sigma,p}$ \thisisrevision{(top)}. The distribution of three-phase cells among this subset is additionally uneven: the minimum and maximum of three-phase cells per process with at least one three-phase cell differ much from the even distribution of the overall sum of three-phase cells within the whole domain \thisisrevision{(bottom)}. The distributions of three-phase cells are shown as a function of the number of time steps for the simulation with $128$ MPI processes; cf.~\protect\reftab{PerformancePLICinFS3D}. The meaning of the $y$-axis is given in the legend.}%
\label{fig:load_imbalance}%
\end{figure}
\refFig{runtime} and \reftab{PerformancePLICinFS3D} gather the results of the performance assessment, where \reffig{load_imbalance} provides additional information. The runtime summed over all MPI-processes shown for the previous and new algorithm in \refFig{runtime} is a direct measure for the costs of the whole simulation, which are reduced by around 20\% with the new sequential PLIC algorithm through a reduction of the positioning time by one order of magnitude. The speed-up of 89-95\% obtained for the reconstruction alone leads to an overall speed-up of 18-21\% in the assessed cases. The main implications can be summarized as follows:
As can be seen from the top row in \reffig{runtime}, the performance gain for the sequential PLIC positioning in three-phase cells is substantial in all cases, persisting throughout the entire simulation, which supports the findings of \reffig{experiment_result}. %

Despite the fact that the three-phase PLIC positioning accounts only for a small fraction of the overall computation time (which decreases with spatial resolution), the speed-up in the overall computation time is above 17\%. The rationale behind this effect can be cast as follows: %
\refFig{mpiprocs_vis} shows that all of the three-phase cells are located in the vicinity of the solid sphere. In the simulation underlying \reffig{load_imbalance}, only between 4 and 12 out of 128 processes contain three-phase cells. This uneven distribution of load forces the remaining MPI processes to stay idle. As can be seen from \reffig{experiment_result}, the proposed algorithm reduces the number of iterations for those cells by a factor of about $20$. Due to the leveraging effect of the load imbalance, this translates to significant speed-up in the overall computation time. %

%
%

%
%
\section{Summary and Conclusions}%
The efficient sequential positioning of two PLIC planes in three-phase cells enclosing given volumes was derived and tested. %
The efficient sequential PLIC positioning requires the cell's polyhedron geometry, the volume fractions of the first and second phase, the material ordering and the two normals of the two planes representing the interfaces in a three-phase cell as the input. %
The polyhedron cell is assumed to not intersect itself, but the algorithm allows the polyhedron to be non-convex. %
The primary phases' interface representation is assumed to have no kinks. %
The outputs of our algorithm are the two signed distances of the planes truncating the correct volumes from the polyhedron cell. %

Recursively employing the \textsc{Gaussian} divergence theorem twice, combined with an elaborate choice of the origin for the volume computation of the truncated volumes of a twice truncated arbitrary polyhedron, the volume computation reduces to a summation of quantities associated to the faces of the original polyhedral cell. %
The new method focusses on the non-degenerate three-phase configurations, where the PLIC planes are not (anti-)parallel, but also shows an efficient solution for the (anti-)parallel degenerate cases based on an existing algorithm by \citet{JCP_2021_fbip} to provide a solution for any three-phase configuration. %
The reference point for the efficient volume computation to solve the non-degenerate cases is chosen on the intersection line (inside or outside the convex hull of the polyhedron) of the two planes truncating the polyhedral cell. This choice eliminates the necessity to establish intersection-dependent connectivity at runtime, as the intersecting planes do not contribute to the volume computation anymore. All computations of the volumes are based on the (un)truncated edges of the known original polyhedron cell. Thus, no costly reestablishing of the connectivity of the residual polyhedron is required after the first truncation. %

Numerical experiments in section~\ref{sec:numerical_experiments} show that, for all instances under consideration, the algorithm admits high efficiency: on average, the secondary plane can be positioned with only 1--2 truncations. %
The performance is almost invariant with respect to the topological configuration (\textit{triple}/\textit{non-wetted}/\textit{fully wetted}), indicating the suitability of the present algorithm for parallel computations with a domain decomposition approach; cf.~\refnote{application_parallel}. %
Our findings highlight the superiority in terms of performance of divergence-based volume computation, which also allows to easily obtain and exploit the derivatives of the parametrized volume. %
In comparison to the decomposition-based bisection approach employed as the reference implementation, the number of truncations required for the positioning was substantially reduced by about one order of magnitude for \textit{triple} configurations. %
The embedded performance assessment in a parallelized CFD solver with a domain decomposition for the parallelization shows that the benefit in a parallelized application exceeds the sole reduction of the computation time by far for a cuboidal domain decomposition like in FS3D: for the configuration considered in this study, a speed-up in the positioning of around 90\% translates to a reduction of the overall computation time by around 20\% despite the fraction of effective runtime spent in the reconstruction is a much smaller percentage. Due to the local concentration of three-phase cells around the primary phase, a load imbalance is induced through the domain decomposition. The load imbalance is thus reduced along with the runtime of the reconstruction, eliminating idle time along with the (effective) reconstruction runtime in the parallelized application. Thus the induced load imbalance becomes insignificant and costly strategies like load rebalancing are not required. %

Beyond the application shown in the present work, the approach does not share the limitation to cuboids as the algorithm does not assume a polyhedron shape other than not self-intersecting. %
The range of applicability of the method covers arbitrary polyhedra with planar faces. %
A reference implementation of the efficient sequential PLIC reconstruction in three-phase cells is available under 
\href{https://doi.org/10.18419/darus-2488}{https://doi.org/10.18419/darus-2488}. The cuboidal cell in the reference implementation can be replaced by any non-self-intersecting polyhedron represented by a list of faces. %
As the performance gain for the sequential PLIC positioning is already large for the whole \textsc{Cartesian} application FS3D, the benefit of a fast, decomposition-free approach for the sequential PLIC positioning without connectivity computation after the first truncation is expected to significantly speed-up three-phase geometric VoF simulations on generalized unstructured polyhedral meshes. %

Concluding, the exploitation of the shared geometry of the original and the truncated residual polyhedron, employing recursive application of the \textsc{Gaussian} divergence theorem, combined with an efficient root-finding leads to the high performance observed for the new algorithm. The application of the \textsc{Gaussian} divergence theorem reduces the computation of truncated volumes to the computation of appropriate segments lying on the edges of the original polyhedron. In the example application, around 20\% of the simulation time can be saved in representative test cases along with the costs and the power consumption of the hardware employed for the multi-phase flow simulations. %

\FloatBarrier

%
%
\begin{center}%
\textsc{Acknowledgment}\\[2ex]%
The authors gratefully acknowledge financial support provided by the German Research Foundation (DFG) within the scope of \href{www.sfbtrr75.de}{SFB-TRR 75 (project number 84292822)}. %
This work was also partly funded by the Deutsche Forschungsgemeinschaft (DFG, German Research Foundation) under the grant ID 265191195 within the \href{https://www.sfb1194.tu-darmstadt.de}{CRC 1194} (Dieter Bothe) and under \href{https://www.simtech.uni-stuttgart.de/exc/research/pn/pn1/}{Germany’s Excellence Strategy - EXC 2075 – 390740016} (Johanna Potyka and Kathrin Schulte).
The simulations in subsection~\ref{Subsec:PerfInFS3D} were conducted on the supercomputer \href{https://www.hlrs.de/solutions/systems/hpe-apollo-hawk}{HPE Apollo (Hawk)} at the \href{https://www.hlrs.de}{High-Performance Computing Center Stuttgart (HLRS)} under the grant no.~FS3D/11142. The authors kindly acknowledge the granted resources and support. 
Furthermore, the authors would like to thank Moritz Heinemann from the \href{https://www.visus.uni-stuttgart.de/}{Visualization Research Center at the University of Stuttgart (VISUS)} for providing the three-phase PLIC \href{https://github.com/UniStuttgart-VISUS/tpf}{visualisation plugin} underlying \reffig{mpiprocs_vis}. Some of the figures and flowcharts in this manuscript were produced using the versatile and powerful library \href{https://ctan.org/topic/pstricks}{\texttt{pstricks}}. For further details and a collection of examples, the reader is referred to the book of \citet{pstricks_2008}. \\[12pt]%

\href{https://www.elsevier.com/authors/policies-and-guidelines/credit-author-statement}{\textsc{CRediT statement}}\\[2ex]%
\textbf{Johannes Kromer}: conceptualization, methodology, software, validation, investigation, data curation, visualisation, writing--original draft preparation, writing--reviewing and editing\\ %
\textbf{Johanna Potyka}: conceptualization, methodology, software, validation, investigation, data curation, visualisation, writing--original draft preparation, writing--reviewing and editing\\ %
\textbf{Kathrin Schulte}: supervision, writing--reviewing and editing, funding acquisition, project administration\\%
\textbf{Dieter Bothe}: conceptualization, methodology, investigation, writing--reviewing and editing, funding acquisition, project administration%
\end{center}%
%
%
\begin{appendix}%
\section{Implicit bracketing}\label{app:implicit_bracketing}%
\begin{figure*}[!tbhp]%
\null\hfill%
\includegraphics[page=1]{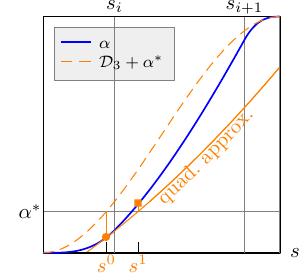}
\hfill%
\includegraphics[page=2]{\thisfigurepath implicit_bracketing}
\hfill%
\includegraphics[page=3]{\thisfigurepath implicit_bracketing}
\hfill\null%
\caption{Implicit bracketing and locally quadratic approximation of the volume fraction $\polyvof\fof{\signdist}$, where $\mathcal{D}_3(\polyvof)$ is a cubic spline with derivative 0 at the boundaries, i.e.\ $\mathcal{D}_3^\prime(s_\pm) =0$.}%
\label{fig:implicit_bracketing_illustration}%
\end{figure*}
\begin{center}
\textit{This section heavily draws from a previous work of the first and last author; \cite[sections~2.3 and 2.4]{JCP_2021_fbip}.}%
\end{center}
The volume fraction $\polyvof\fof{\signdist}$ is an increasing cubic polynomial, denoted by $\mathcal{S}_i\fof{\signdist}$, within a bracket $\mathcal{B}_i\defeq[\signdist_i,\signdist_{i+1}]$. The polynomial reads %
\begin{equation}
\begin{aligned}
\mathcal{S}_i\fof{z;\signdist^n} = &
\frac{\polyvof^{\prime\prime\prime}\fof{\signdist^n}}{6} (z-\signdist^n)^3  + \frac{\polyvof^{\prime\prime}\fof{\signdist^n}}{2} (z-\signdist^n)^2 \\
  & + \polyvof^{\prime}\fof{\signdist^n}(z-\signdist^n)  + \polyvof\fof{\signdist^n} 
\end{aligned}
\label{eqn:third_order_approximation}%
\end{equation}%
for any given $\signdist^n\in\mathcal{B}_i$. %
Hence, the truncation of the polyhedron $\polyhedron*$ at any $\signdist^n$ (implicitly) provides the full information of $\polyvof\fof{\signdist}$ within the containing bracket $\mathcal{B}_i$. Exploiting that $\polyvof_i=\mathcal{S}_i\fof{\signdist_i;\signdist^n}$ and $\polyvof_{i+1}=\mathcal{S}_i\fof{\signdist_{i+1};\signdist^n}$ suggests the following strategy: %
\begin{enumerate}%
\item if the current iteration $\signdist^n$ is \textit{not} located in the target bracket $\refbracket$ ($\refvof<\polyvof_i$ or $\polyvof_{i+1}<\refvof$), the next iteration is obtained from locally quadratic approximation (\reffig{implicit_bracketing_illustration}, left). %
\item if the current iteration $\signdist^n$ is located in the target bracket $\refbracket$ ($\polyvof_i\leq\refvof\leq\polyvof_{i+1}$), the sought $\signdistref$ corresponds to the root of $\mathcal{S}_i-\refvof$, requiring no further truncation (\reffig{implicit_bracketing_illustration}, center). %
\end{enumerate}%

The initial guess is obtained from a global cubic spline interpolation via %
\begin{equation}
\begin{aligned}
s^0\fof{\refvof}\defeq & \projvertposmin+\brackets{\projvertposmax-\projvertposmin}\Bigl(\frac{1}{2}-\\
 & \cos\brackets*{\frac{\arccos\fof{2\refvof-1}-2\pi}{3}}\Bigr)
\end{aligned}
\end{equation}
with $\projvertposmin=\min\fof{\hat{\mathcal{S}}}$ and $\projvertposmax=\max\fof{\hat{\mathcal{S}}}$; cf.~\refeqn{vertex_normal_projections}. %
\refFig{implicit_bracketing_illustration} illustrates the components of the strategy outlined above: with $\signdist^0\not\in\mathcal{B}_i$ (left), a quadratic approximation yields $\signdist^{1}$, which lies within the target bracket $\refbracket=\mathcal{B}_i$ (center). %
The rightmost configuration already starts with $\signdist^0\in\refbracket$, such that the spline interpolation directly yields the sought position $\signdistref$, implying that only a single truncation is required. %
\section{An explicit inverse of the cuboid volume function}\label{app:cube_explicit_volume_inverse}%
Note that any cuboid with edge lengths $\set{\Delta x_i}_{i=1}^{3}$ can be transformed to the reference cube with $\Delta x_i\equiv1$. Furthermore, one can always find a rotation such that the components of the normal $\plicnormal$ are positive and arranged in descending order. %
Let $n_i\defeq\iprod{\plicnormal}{\ve_i}$ for ease of notation. %
Assuming that the second and third component of $\plicnormal$ are zero ($n_2=n_3=0$), the positioning problem becomes degenerate and one trivially obtains $\signdistref=\refvof\Delta x$. The normal can be parametrized by the polar angle $\varphi$ with $\tan\varphi=\frac{n_2}{n_1}$ for the special case of two non-zero components ($n_3=0$), and one obtains %
\begin{equation}
\begin{aligned}
\signdistref=\Delta x%
\begin{cases}
\sqrt{\refvof\sin2\varphi} \text{ for } \refvof\leq\frac{\tan\varphi}{2} \text{,}\\[10pt]
\sin\varphi+ \cos\varphi(\refvof- \frac{\tan\varphi}{2}) \\
\quad \text{ for } \frac{\tan\varphi}{2}<\refvof<1-\frac{\tan\varphi}{2}\text{,}\\[10pt]
\cos\varphi+ \sin\varphi-\sqrt{\brackets*{1-\refvof}\sin2\varphi} \\
\quad \text{ for } \refvof\geq1-\frac{\tan\varphi}{2}\text{,}%
\end{cases}
\end{aligned}
\end{equation}
where, due to the symmetry, the second case degenerates for $\varphi=\nicefrac{\pi}{4}$. %
If all components of $\plicnormal$ are non-zero, the signed distance reads %
\begin{equation}
\signdistref = \mathcal{S}_k\fof{\refvof}\quad \text{for} \quad\refvof\in(\hat{\polyvof}_{k} , \hat{\polyvof}_{k+1}]
\label{eqn:calclstar_vof}
\end{equation}
with%
\begin{equation}
\begin{aligned}
\vec{\hat{\polyvof}}=%
\begin{cases}
\Bigl[0, &  \frac{n_3^2}{6n_1n_2},   \frac{n_2^3-(n_2-n_3)^3}{6n_1n_2n_3} \frac{n_2+n_3}{2n_1}, \frac{1}{2}\Bigr] \\ 
 & \text{ for }\quad n_1\geq n_2+n_3,\\[10pt]%
\Bigl[0, &  \frac{n_3^2}{6n_1n_2},  \frac{n_2^3-\brackets{n_2-n_3}^3}{6n_1n_2n_3}, \frac{n_1^3-\brackets{n_1-n_2}^3-\brackets{n_1-n_3}^3}{6n_1n_2n_3}, \frac{1}{2}\Bigr] \\
& \text{ for }\quad n_1<n_2+n_3,%
\end{cases}
\end{aligned}\label{eqn:calclstar_vof_limits}%
\end{equation}
and %
\begin{align*}
\mathcal{S}_1\fof{\refvof}= \sqrt[3]{6\refvof n_1 n_2 n_3}\text{,}%
\end{align*}
\begin{align*}
\mathcal{S}_2\fof{\refvof}=\frac{n_3}{2}+\sqrt{2 \refvof n_1 n_2-\frac{n_3^2}{12}}\text{,}%
\end{align*}
\begin{align*}
\mathcal{S}_3\fof{\refvof}=&n_2+n_3-\sqrt{8n_2n_3}\cos\Biggl(\frac{1}{3}\arccos\Bigl(\frac{3}{8}\sqrt{\frac{2}{n_2n_3}}\brackets{n_2 \\
& + n_3-2n_1\refvof}\Bigr)+\frac{\pi}{3}\Biggr),\\%
\end{align*}
\begin{align*}
\mathcal{S}_4\fof{\refvof}&=%
\begin{cases}
\frac{2\refvof n_1+n_2+n_3}{2} & \quad \text{ for } n_1\geq n_2+n_3\\[10pt]%
\frac{n_1+ n_2+n_3}{2} + & 2\sqrt{\frac{-p}{3}}\cos(\frac{1}{3}\arccos\brackets*{\frac{3q}{2p}\sqrt{\frac{-3}{p}}} \\ 
 & -\frac{2\pi}{3}) \quad \text{ for } n_1<n_2+n_3%
\end{cases}%
\end{align*}
with %
\begin{equation*}
p = \frac{3}{4}(2 n_1^2-(n_1 + n_2 + n_3)^2 + 2 n_2^2 + 2 n_3^2)
\end{equation*}
and
\begin{equation*}
q =\frac{3 n_1 n_2 n_3}{2}(2 \refvof-1) \text{.}
\end{equation*}
Note that the point-symmetry of the cube volume fraction $\refvof$ can be exploited to extend \refeqn{calclstar_vof} to $\refvof>\frac{1}{2}$; cf.~\reffig{calclstar_cuboid_illustration} for an illustration. A similar set of formulae was derived in \citet{JCP_2000_arcl} and \citet{arxiv_2020_astt}. %
The volume's inversion is included here, as it is essential for the previous three-phase PLIC positioning algorithm employed in FS3D described in appendix~\ref{app:decomposition_approach} and used in this paper for comparison.
\begin{figure*}%
\null\hfill%
\includegraphics[page=1]{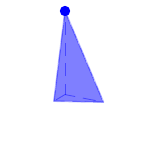}%
\hfill
\includegraphics[page=2]{\thisfigurepath calclstar_illustration_3D}%
\hfill
\includegraphics[page=3]{\thisfigurepath calclstar_illustration_3D}%
\hfill
\includegraphics[page=4]{\thisfigurepath calclstar_illustration_3D}%
\hfill
\includegraphics[page=5]{\thisfigurepath calclstar_illustration_3D}%
\hfill
\includegraphics[page=6]{\thisfigurepath calclstar_illustration_3D}%
\hfill
\null%
\\%
\null\hfill%
\includegraphics[page=7]{\thisfigurepath calclstar_illustration_3D}%
\hfill%
\includegraphics[page=8]{\thisfigurepath calclstar_illustration_3D}%
\hfill%
\includegraphics[page=9]{\thisfigurepath calclstar_illustration_3D}%
\hfill%
\includegraphics[page=10]{\thisfigurepath calclstar_illustration_3D}%
\hfill%
\includegraphics[page=11]{\thisfigurepath calclstar_illustration_3D}%
\hfill%
\includegraphics[page=12]{\thisfigurepath calclstar_illustration_3D}%
\\%
\null\hfill%
\includegraphics[page=1]{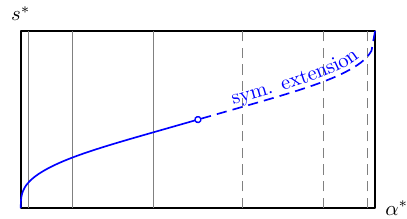}%
\hfill%
\includegraphics[page=2]{\thisfigurepath calclstar_vof_3D}%
\hfill\null%
\caption{Intersection topology for cube with normal $\plicnormal=\frac{[4,2,1]\transpose}{\sqrt{21}}$ ($n_1>n_2+n_3$, blue) and $\plicnormal=\frac{[4,3,2]\transpose}{\sqrt{29}}$ ($n_1<n_2+n_3$, red) at $\hat{\polyvof}_k$ (vertical lines) from \protect\refeqn{calclstar_vof_limits} with the respectively intersected vertices ($\bullet$).}%
\label{fig:calclstar_cuboid_illustration}%
\end{figure*}%
\end{appendix}%
\section{A decomposition-based accelerated bisection approach for cuboids}\label{app:decomposition_approach}%
\begin{figure*}[htbp]
\null\hfill%
\includegraphics{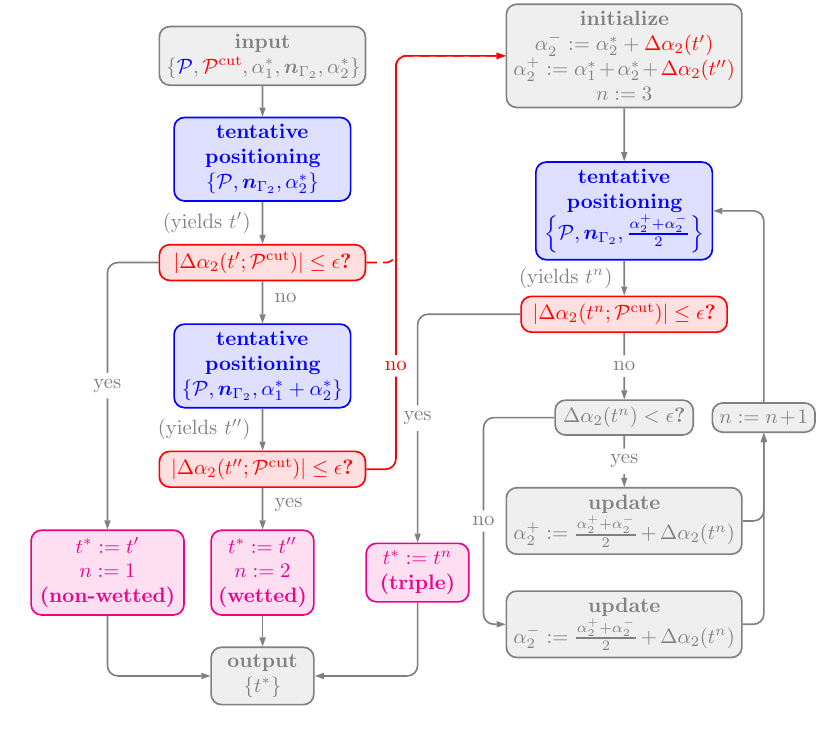}%
\hfill\null%
\caption{Flowchart of the accelerated decomposition-based bisection approach of \citet{Potyka_2018} with $\Delta\polyvof_2\fof{\secondarysigndist} = \Delta\polyvof_2\fof{\secondarysigndist,\cutpolyhedron} \protect\defeq\secondaryrefvof-\polyvof_2\fof{\secondarysigndist;\cutpolyhedron}$ and $\epsilon=\num{e-12}$; cf.~\protect\reffig{decomposition_vof_error} for an illustration. $\polyvof_2^+$ and $\polyvof_2^-$ reffer to the whole cell $\protect\polyhedron*$, thus the tentative positioning is performed with the two-phase algorithm, cf.~\refapp{cube_explicit_volume_inverse}.}%
\label{fig:flowchart_decomposition}
\end{figure*}

\begin{figure*}[htbp]
\null\hfill%
\includegraphics[page=1]{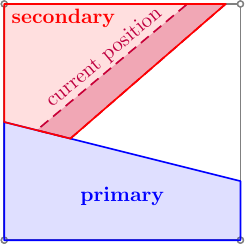}%
\hfill%
\includegraphics[page=3]{\thisfigurepath literature_decomposition}%
\hfill%
\includegraphics[page=2]{\thisfigurepath literature_decomposition}%
\hfill%
\includegraphics[page=4]{\thisfigurepath literature_decomposition}%
\hfill\null%
\caption{Update strategy for decomposition-based approach of \citet{Potyka_2018}: effectively, the updated position $\secondarysigndist^{n+1}$ is obtained by re-distributing the deviation $\Delta\polyvof_2\fof{\secondarysigndist^n;\cutpolyhedron}$. Note that the domain overlap and the deviation of the next iteration (crosshatch) are coextensive, i.e., they admit the same volume. In other words, as the iterative position $\secondarysigndist^n$ approaches the true position $\secondarysigndistref$, the associated deviation vanishes, i.e., $\Delta\polyvof_2\fof{\secondarysigndist^n;\cutpolyhedron}\to0$ as $\secondarysigndist^n\to\secondarysigndistref$. %
The error $\Delta \polyvof_2\fof{\secondarysigndist^n;\cutpolyhedron}$ at the same time is a lower bound for the required change of the volume fraction enclosed inside the whole cell $\polyvof_2\fof{\secondarysigndist^n;\protect\polyhedron*}$.} \label{fig:decomposition_vof_error}%
\end{figure*}
The three-phase PLIC positioning algorithm formerly employed in FS3D and used for comparison in this work was described in a Master's Thesis in German \cite{Potyka_2018}, thus a short version of it's description is provided in the following: %
The idea of this decomposition-based sequential positioning in three-phase cells is, that the positioning for two-phase cells is reused. This is feasible, if the positioning in two-phase cells is fast. In the case of a \textsc{Cartesian} mesh the two-phase PLIC position is retrieved from an explicit formula, cf.~\refapp{cube_explicit_volume_inverse}. Therefore the requirement of a fast algorithm for two-phases is fulfilled. Additionally, the first PLIC plane is positioned with the two-phase algorithm shown in \refapp{cube_explicit_volume_inverse} which is a function of the volume fraction $\polyvof_1$, the normal $\primaryplicnormal$ and the grid cell's dimensions $\Delta x$. An iterative algorithm, which reuses the two-phase positioning, is utilized to position the secondary PLIC plane. The initial minimum and maximum value for $\secondarysigndist$ are defined by the two extrema of the interface's orientations where no triple line occurs:
\begin{enumerate}
\item The two interfaces do not intersect or cover each other inside the cell, they are independent, thus $\secondarysigndist'=\mathcal{S}\fof{\polyvof_2^*}$ and $\polyvof_{2}^-=\polyvof_2\fof{\secondarysigndist';\polyhedron*}$
\item The secondary plane fully wets the primary plane, thus $\secondarysigndist''=\mathcal{S}\fof{\polyvof_1^* + \polyvof_2^*}$ and $\polyvof_{2}^+=\polyvof_2\fof{\secondarysigndist'';\polyhedron*}$.
\end{enumerate}
\begin{figure*}
\null\hfill
\includegraphics[page=1]{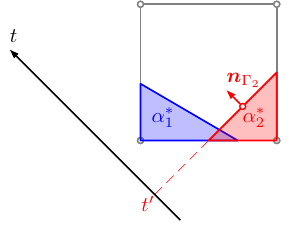}%
\hfill
\includegraphics[page=2]{\thisfigurepath face_bracket_decomposition}%
\hfill \null
\caption{The initial bracketing values for the accelerated bisection are found with the assumption of a \textit{non-wetted} configuration for the lower limit $\secondarysigndist^\prime$ and a \textit{wetted} assumption for upper limit $\secondarysigndist^{\prime\prime}$. If one of the assumed topologies is present, no further iteration is necessary. The secondary PLIC plane's position $\secondarysigndistref$ of a \textit{triple} configuration is bound to lie between those extrema. Thus, $\primaryrefvof$ determines the size of the initial bracket.} %
\label{fig:face_bracket_decomposition}%
\end{figure*}
The initial bracketing values for $\polyvof_2\fof{\secondarysigndist;\polyhedron*}$ and their corresponding positions $\secondarysigndist$ are depicted in \reffig{face_bracket_decomposition}. Both limits can be retrieved from the explicit two-phase PLIC positioning method in cuboids applied to (a sum of) the target volumes, cf.~\refapp{cube_explicit_volume_inverse}. For each tentative position $\secondarysigndist^n$ from the volume fraction $\polyvof_2\fof{\secondarysigndist^n;\polyhedron*}$ enclosed inside the whole cell $\polyhedron*$ the volume fraction $\polyvof_2\fof{\secondarysigndist^n;\cutpolyhedron}$ inside the residual polyhedron $\cutpolyhedron$ is computed. The error in the enclosed volume
\begin{equation}
\Delta \polyvof_2\fof{\secondarysigndist^n,\cutpolyhedron} = \polyvof_2^* - \polyvof_2\fof{\secondarysigndist^n,\cutpolyhedron}
\label{eqn:ferror}
\end{equation}
of the present secondary position is compared to the true volume fraction $\polyvof_2^*$. If either the secondary position at the minimum $\polyvof_2^-$ or maximum $\polyvof_2^+$ volume fraction inside the whole cell $\polyhedron*$ already yield an error below a given tolerance, a further iteration of the second PLIC plane's position is not necessary, i.e., the liquid does not wet or fully wet the solid as no intersection volume is computed. If a contact line exists, the required tolerance is not reached with those two extreme positions. The calculated error is reused to refine the initial values enclosed by the PLIC plane inside the whole cell, as the \refeqn{ferror} is also a lower limit for the required adjustment in the volume enclosed below the secondary plane inside the whole cell $\polyvof_2\fof{\secondarysigndist^n;\polyhedron*}$ for the next iteration, cf.~\reffig{decomposition_vof_error}. Depending on the sign of $\Delta \polyvof_2\fof{\secondarysigndist^n;\cutpolyhedron}$, either the upper bound $\polyvof_2^+$ or the lower bound $\polyvof_2^-$ is adjusted with $\Delta \polyvof_2\fof{\secondarysigndist^n;\cutpolyhedron}$. %
Both volume fractions are corresponding to the extended PLIC planes, the enclosed volume fraction in $\polyhedron*$. Therefore for this updated volume fraction \refapp{cube_explicit_volume_inverse} is applicable to calculate the next update of $\secondarysigndist^n$. A faster convergence is reached by combining a bisection based update with the calculated error. In each iteration, the volume fraction's error \refeqn{ferror} is evaluated and $\polyvof_{2}^-$ or $\polyvof_{2}^+$ updated depending on the sign of $\Delta \polyvof_2\fof{\secondarysigndist^n;\polyhedron*}$ like described above until the given tolerance is fulfilled in addition to employing a minimisation with bisection. The calculation of the error $\Delta \polyvof_2\fof{\secondarysigndist^n;\cutpolyhedron}$ is depicted in \reffig{decomposition_vof_error}. The flowchart in \reffig{flowchart_decomposition} illustrates the described procedure to determine the planes' positions.\\%
For robust convergence a few special cases are considered: %
\begin{itemize}
\item A tolerance is needed to ensure that $\polyvof_2\fof{\secondarysigndist^n;\polyhedron*} + \Delta \polyvof_2\fof{\secondarysigndist^n;\cutpolyhedron}$ does not overshoot or undershoot the minimum and maximum boundaries, if the new volume fraction corresponds to $\secondarysigndistref$ except for rounding errors during the computation.
\item  If the polyhedron's dimensions become as small that there is up to a tolerance no difference between the minimum or maximum vertex point of the cell and an intersection point, $\polyvof_2\fof{\secondarysigndist^n;\polyhedron*} = 0$ is set at the minimum or $\polyvof_2\fof{\secondarysigndist^n;\polyhedron*} = 1$ at the maximum for the initial bracketing values. %
\end{itemize}
The volume error computation consists of an intersection of the residual polyhedron $\cutpolyhedron$ with $\secondarysigndist^n$ and the volume computation itself. For the computation of the intersection, all vertices are marked, if they are above, below or on the plane. All points which are below the plane and on the plane belong to the new polyhedron, but also new points have to be computed which are intersection points of the plane with the segments where one point is above and one below the plane. The connectivity of the original and intersected polyhedron are computed for the compuation of a new point. This requires the information of common faces, thus multiple comparison operations are necessary. The volume computation of the polyhedron exploits, that the polyhedrons resulting from sequential cuts of a cuboid with one or two planes is convex. From an origin inside the polyhedron the volumes of tetrahedra are summed up. For each tetrahedron the connectivity has to be computed again by multiple comparisons, if the corner points have common faces. This requires significant computational effort already for cuboids with two intersections resulting in convex polyhedrons.
\section{FS3D Simulation Setup}\label{app:FS3Dcollisionsetup}
\begin{figure*}[bt!]%
\null\hfill
\begin{minipage}[b]{0.6\textwidth}
\centering
\includegraphics{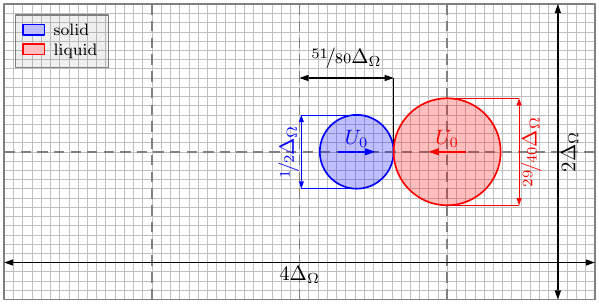}
\end{minipage}
\hfill%
\begin{minipage}[b]{0.3\textwidth}
\centering
 \begin{tabular}[b]{c|c|c}
 $\Delta_\Omega$ & $2/5$ & $\mathrm{cm}$ \\
 $U_0$           & $90$  & $\mathrm{cm/s}$ \\
 \hline
 $\rho_s$           & $2.55$  & $\mathrm{g/cm^3}$ \\
 $\rho_l$           & $1$  & $\mathrm{g/cm^3}$ \\
 $\mu_l$            & $0.01$  & $\mathrm{g/(cm s)}$ \\
 $\rho_g$           & $1.19 \cdot 10^{-3}$  & $\mathrm{g/cm^3}$ \\
 $\mu_g$            & $1.824 \cdot 10^{-4}$   & $\mathrm{g/(cm s)}$ \\
 \hline
 $\sigma_{lg}$      &  $72.0$ & $\mathrm{g/s^2}$ \\
 $S$                &  $53.5$ & $\mathrm{g/s^2}$ \\
 \end{tabular}
 \vspace{\baselineskip}
 \end{minipage}
\hfill\null%
\caption{Setup of the test case (center section): collision of a super-hydrophilic solid sphere and a liquid droplet in the domain $\Omega={[-2\Delta_\Omega,2\Delta_\Omega]\times[-\Delta_\Omega,\Delta_\Omega]^2}$. To reduce computation  time, at $t=0$ the particles initially touch, where the offset ensures that the ligament stretches across the entire domain. The dashed lines indicate the cuboidal domain decomposition ($2P\times P\times P$ processes, shown for $P=2$) underlying the parallelization. The table provides the initialized values of the half domain width $\Delta_\Omega$ (c.f., sketch on the left), the velocity $U_0$ of both the particle and the droplet as well as the employed liquid's (index $l$), solid's (index $s$) and gas's (index $g$) properties: the densities $\rho$, viscosities $\mu$ and surface tension $\sigma_{lg}$ are given. The spreading parameter $S >0$ shows that a fully wetting surface-liquid combination is simulated.
}%
\label{fig:FS3Dsetup}%
\end{figure*}
The setup shown in \reffig{FS3Dsetup} of an artificial collision of a solid sphere (subscript s) and a liquid droplet (subscript l) in ambient gas (subscript g) is used for the prototypical test cases in this work. The simulation setup is chosen such, that a qualitatively realistic volume distribution for the collision with a hydrophilic sphere is reached. The simulations are not validated against experiments and do not necessarily show the true physical process. The focus of this test setup was to produce a representative volume distribution and spectrum of three-phase cells in order to test the new reconstruction algorithm.

\FloatBarrier
%
%
\renewcommand{\bibname}{References}
\bibliographystyle{abbrvnat}
\bibliography{literature_full}%

\end{document}
\endinput